\documentclass[12pt,oneside,english]{amsart}
\usepackage{tcolorbox}
\usepackage{graphicx}
\usepackage{caption}

\def\ov{\mathrm{overlap}}
\def\iv{\mathrm{agree}}
\newcommand{\Z}{\mathbb{Z}}
\newcommand{\B}{\euscript{B}}
\usepackage{amsmath}
\usepackage{amssymb}
\usepackage{amsthm}
\usepackage{hyperref}
\usepackage{enumerate}
\usepackage{graphicx}
\usepackage{color}
\usepackage{xcolor}
\usepackage{mathtools}
\usepackage{float}
\usepackage{tikz}
\usetikzlibrary{arrows}
\usepackage{bm}
\usepackage[title]{appendix}
\usepackage{mathrsfs}
\DeclareMathAlphabet\euscript{U}{eus}{m}{n}
\usepackage{hyperref}
\hypersetup{
	colorlinks = true, 
	urlcolor = blue, 
	linkcolor = red, 
	citecolor = blue 
}

\newtheorem{thm}{Theorem}[section]
\newtheorem{co}[thm]{Corollary}
\newtheorem{lem}[thm]{Lemma}
\newtheorem{ques}[thm]{Question}
\newtheorem{prop}[thm]{Proposition}
\newtheorem{conj}[thm]{Conjecture}
\newtheorem{assumption}[thm]{Assumption}

\newtheorem{pr}[thm]{Proposition}

\newtheorem{definition}[thm]{Definition}
\newenvironment{de}{\begin{definition}\rm}{\end{definition}}
\newtheorem{example}[thm]{Example}
\newenvironment{exmp}{\begin{example}\rm}{\end{example}}
\newtheorem{remark}[thm]{Remark}
\newenvironment{rem}{\begin{remark}\rm}{\end{remark}}
\newtheorem{algorithm}[thm]{Algorithm}

\newtheorem{observation}[thm]{Observation}

\newtheorem{claim}[thm]{Claim}

\newtheorem{fact}[thm]{Fact}

\newcommand\blfootnote[1]{%
	\begingroup
	\renewcommand\thefootnote{}\footnote{#1}%
	\addtocounter{footnote}{-1}%
	\endgroup
}

\newcommand{\vertiii}[1]{{\left\vert\kern-0.25ex\left\vert\kern-0.25ex\left\vert #1
    \right\vert\kern-0.25ex\right\vert\kern-0.25ex\right\vert}}
\title{Shifts of Finite Type Obtained by Forbidding a Single Pattern}

\author{
	Nishant Chandgotia$^1$
	\and
	Brian Marcus$^2$
	\and
	Jacob Richey$^3$
	\and
	Chengyu Wu$^4$
}

\address{Tata Institute of Fundamental Research - Centre for Applicable Mathematics, Bengaluru}
\email{nishant.chandgotia@gmail.com}
\address{The University of British Columbia, Vancouver}
\email{marcus@math.ubc.ca}
\address{The Alfr\'ed R\'enyi Institute of Mathematics, Budapest}
\email{jfrichey001@gmail.com}
\address{The University of British Columbia, Vancouver}
\email{chengyuw@connect.hku.hk}

\begin{document}\thispagestyle{empty}
\begin{abstract}
Given a finite word $w$, Guibas and Odlyzko (\emph{J. Combin. Theory Ser. A, 30, 1981, 183-208}) showed that the autocorrelation
polynomial $\phi_w(t)$ of $w$, which records the set of self-overlaps of $w$, explicitly
determines for each $n$, the number $|\euscript{B}_n(w)|$  of words of length $n$ that avoid $w$.
We consider this and related problems from the viewpoint of symbolic dynamics, focusing on the setting of $X_{\{w\}}$, the space of all bi-infinite sequences that avoid $w$.  We first summarize and elaborate upon (\emph{J. Combin. Theory Ser. A, 30, 1981, 183-208}) and other work to show that the sequence $|\euscript{B}_n(w)|$  is equivalent to several invariants of $X_{\{w\}}$.
We then give a finite-state labeled graphical representation $L_w$ of $X_{\{w\}}$ and show that $w$ can be recovered from the graph isomorphism class of the unlabeled version of $L_w$. Using $L_w$, we apply ideas from probability and Perron-Frobenius theory to obtain results comparing features of $X_{\{w\}}$ for different $w$.
 Next, we give partial results on the problem of classifying the spaces $X_{\{w\}}$ up to conjugacy.  Finally, we extend some of our results to spaces of multi-dimensional arrays that avoid a given finite pattern.
\end{abstract}
\maketitle

%

\section{Introduction}\label{section:intro}
\bigskip
\blfootnote{$^1$Tata Institute of Fundamental Research-Centre for Applicable Mathematics, $^2$ The University of British Columbia, $^3$The Alfr\'ed R\'enyi Institute of Mathematics, $^4$ The University of British Columbia}
For a positive integer $q$ and a finite word $w = w_1 \cdots w_{k}$ over the alphabet
$[q]:= \{0, \ldots, q-1\}$, consider the set $\euscript{B}_n(w)$ of words over $[q]$ of length $n$ which do {\em not} contain $w$, i.e.,  in which the  word $w$ is forbidden as a contiguous subword. A natural problem is to determine how the sequence of numbers $|\euscript{B}_n(w)|$ depends on the word $w$.  A  complete answer to this problem was given over forty years ago in seminal work by Guibas and Odlyzko~\cite{GO81}: the starting point is an explicit formula for the generating function
$G_w(t) = \sum_{n=0}^\infty |\euscript{B}_n(w)|t^n$. Their formula is a simple function of the {\em (auto-) correlation polynomial}, $\phi_u(t)$, whose coefficients record the lengths of self-overlaps of $u$. The formula shows that
\begin{equation}
\label{one}
G_u(t) = G_v(t) ~\mbox{iff}~ \phi_u(t) = \phi_v(t).
\end{equation}
Since the sequence  $|\euscript{B}_n(w)|$ grows at most exponentially in $n$, the correlation polynomial $\phi_w(t)$ completely determines the sequence $|\euscript{B}_n(w)|$.

A natural setting for problems of this type is symbolic dynamics.  For a positive integer $q$,  the {\em full $q$-shift} is the set $[q]^\mathbb{Z}$ of bi-infinite sequences over the alphabet $[q]$. A {\em shift space} over $[q]$ is a subset $X \subset [q]^\mathbb{Z}$ defined by a list $\mathcal{F}$ (finite or infinite) of finite forbidden words.  A prominent example is $X = X_{\{w\}}$, the set of all bi-infinite sequences over $[q]$ that do not contain the (finite) word $w$.   This combines all the $\euscript{B}_n(w)$ above neatly into a single object. A simple example is to take $q=2$ and $w=11$, and then $X_{\{w\}}$ is called the {\em golden mean shift.} More generally one can consider a {\em shift of finite type (SFT)} which is a shift space defined by a finite list of finite forbidden words. Of course, $X_{\{w\}}$ is an SFT defined by a single forbidden word $w$. Another  prominent class of SFTs is defined by finite directed graphs $G$: the {\em edge shift} $X_G$ is the set of bi-infinite edge walks on $G$.

Of particular interest in symbolic dynamics is the {\em entropy} $h(X)$ which is the asymptotic growth rate of the number of allowed finite sequences in elements of the shift space $X$.
Another useful object is the {\em zeta function} $\zeta(X)$, defined in Section~\ref{section:follower}, which is a function of a variable $t$ that encapsulates the number of periodic sequences in the shift space $X$ of all periods. We often write $\zeta(X_{\{w\}})$ as $\zeta_w(t)$. See~\cite{LM21} for more background on symbolic dynamics.

In Section~\ref{section:earlier_1d}, we
piece together results from Guibas and Odlyzko~\cite{GO81}, Lind~\cite{Lin89}, and Eriksson~\cite{Eri97}, to show in Proposition \ref{equality},  that  $\phi_w(t)$, $G_w(t)$,  $h(X_{\{w\}})$, $\zeta(X_{\{w\}})$, as well as a  probabilistic quantity $E(\tau_w)$,  all determine essentially the same information (here $E(\tau_w)$ is the expected time to generate the word $w$ in a process that generates symbols in $[q]$ uniformly and independently). In Proposition \ref{inequality_1d}  we summarize related results involving inequalities in these quantities. Our aim in both propositions is to clarify the various implications with simple proofs: to do so we mostly rely on statements and methods of previous authors, plus a few new ideas.

One of the main results of our paper (Theorem \ref{Theorem: replace language equal} of Section~\ref{section:higherd}) is a far-reaching generalization of the sufficiency (``if'' part) of (\ref{one}) above, using an entirely different method.
Consider a finite subset $S$ of the $d$-dimensional lattice $\Z^d$. A {\em pattern} on $S$ is an element $w \in [q]^S$. We consider an analogue of the self-overlap set, called the {\em self-agreement set} of $w$ on $S$.

Given such a $w$ and another  finite subset $T \subset \Z^d$, one can consider the set $\euscript{B}_T(w)$ of patterns $x \in [q]^T$ which do not contain any translate of $w$. We show that for any given $S$ and $u,v \in [q]^S$, if the self-agreement sets of $u$ and $v$ on $S$ are the same, then for all $T \subset \Z^d$, $|\euscript{B}_T(u)| = |\euscript{B}_T(v)|$.  The special case $d=1$ is the setting described above in terms of words where $S = [0,k-1], T = [0,n-1]$. A much harder question is to compare $|\euscript{B}_T(u)|$ and $|\euscript{B}_T(v)|$ when the patterns $u,v \in [q]^S$ don't necessarily have the same set of overlaps. While we cannot address this completely, we mention some preliminary ideas in Section \ref{subsection: higher d inequalities}.

The proof of this generalization is obtained using a simple inclusion-exclusion argument by replacing occurrences of $u$ by $v$ in patterns on finite subsets. This seems to be related to ideas in Jacquet and Szpankowski \cite[Section 4.3]{jacquet2015analytic}; see also \cite{carrigan2024natural}.
 We also consider periodic configurations and show that whenever $u,v \in [q]^S$ have the same ``self-overlap'' sets,
 then for any finite fundamental domain $D \subset \Z^2$, the number of periodic configurations with fundamental domain $D$ that avoid $u$ is the same as the number that avoid $v$.  These results can be generalized to arbitrary dimension and possibly to the setting of shift spaces over arbitrary countable groups. 

Returning to the case $d=1$, given a word $w$ we consider finite-state directed graphical representations of $X_{\{w\}}$ in Section~\ref{section:follower}.  This point of view allows us to make use of Perron-Frobenius theory.
We give a simple explicit description of an edge-labelled graph $L_w$ whose edge sequences are in one-to-one correspondence with elements of $X_{\{w\}}$.  Here $L_w$ is a simple variation of the well-known follower set graph~\cite[Section 3.2]{LM21} from symbolic dynamics.


We show that the graph isomorphism class of the unlabelled version of $L_w$ completely determines the word $w$, up to a possible permutation of the alphabet (Corollary \ref{recover_word}).
We also give certain graphical conditions on $L_u$ and $L_v$ which imply $h(X_{\{u\}}) \ge h(X_{\{v\}})$.  This shows how comparisons between different $X_{\{w\}}$ can be viewed through a graphical lens, involving Perron-Frobenius theory.
We show that the same graphical conditions imply $E(\tau_u) \ge E(\tau_v)$; in fact, the graphical conditions imply that $\tau_u$ stochastically dominates $\tau_v$.
%

For shift spaces $X$ and $Y$, a {\em sliding block code} $\phi:X \to Y$ is a stationary coding  $\phi(x) = y$, with each symbol $y_n$ determined by a ``local'' rule  $y_n = \Phi(x_{[n-m, n+a]})$, for some fixed $m$, $a$.
 A {\em conjugacy} from one shift space $X$ to another $Y$ is a bijective sliding block code from $X$ onto $Y$, in which case one says that $X$ and $Y$ are {\em conjugate.}
  The main classification problem in symbolic dynamics, still unsolved after many years, is the problem of deciding, given two SFTs $X$ and $Y$, whether $X$ and $Y$ are conjugate.

   It is not hard to see that given finite directed graphs $G$ and $H$, the conjugacies $\phi$ from $X_G$ to $X_H$ with $m=a=0$ for both $\phi$ and $\phi^{-1}$ are precisely those given by graph isomorphisms from $G$ to $H$.  But these are rather trivial conjugacies.
It is also not hard  to see that for any finite word $w$ over  $[q]$, $X_{\{w\}}$ is conjugate to the edge shift $X_{L_w}$, in a natural way.  Thus, any graph isomorphism from $L_u$ to $L_v$ defines a conjugacy from $X_{\{u\}}$ to $X_{\{v\}}$.  Since the graph isomorphism class of $L_u$ determines $u$ (modulo a permutation of the letters) there doesn't exist a conjugacy of this type from $X_{\{u\}}$ to $X_{\{v\}}$ (if $u$ and $v$ are essentially different).

Are there other conjugacies?  Since the entropy $h(X_{\{u\}})$ determines the correlation polynomial $\phi_u(t)$ and entropy is preserved under conjugacy, a necessary condition for $X_{\{u\}}$ to be conjugate to $X_{\{v\}}$ is that $\phi_u(t) = \phi_v(t)$.  In the case where $\phi_u(t) = \phi_v(t) = t^{k-1}$, $u$ and $v$ have no non-trivial self-overlaps, and in Theorem \ref{conj_trivial_self_over} of Section~\ref{section:conjugacy}, we show that $X_{\{u\}}$ and $X_{\{v\}}$ are conjugate (under an irreducibility condition).  So, from the viewpoint of conjugacy of $X_{\{u\}}$, the words with no non-trivial self-overlaps are all the ``same,'' but from the viewpoint of graph isomorphism of $L_u$, these words are all ``different.''

As a simple example, $u=110100$ and $v=111010$ have no non-trivial self-overlaps, and so
$X_{\{u\}}$ and $X_{\{v\}}$ are conjugate.
If $u = 110110$ and $v=011011$ (the reverse of $u$), then $u$ and $v$ have the same non-trivial set of self-overlaps,  but we do not know if $X_{\{u\}}$ and $X_{\{v\}}$ are conjugate.

Instead of starting with a full shift and forbidding a single word, one could start with another simple ({\em ambient}) SFT $X$ and then forbid a single word from $X$. For instance, $X$ could be the SFT defined by forbidding a single word $w$ of length two, $X = X_{\{w\}}$, in which case $(X_{\{w\}})_{\{u\}}$ would be the SFT obtained by forbidding $w$ and $u$. A few results in this framework, with $w = 11$, are sprinkled throughout the paper (see Sections~\ref{fol_mem1} and \ref{conj_subsec2}). In fact, this is the framework considered by Lind~\cite{Lin89}, where $w$ is an arbitrary word of length 2.


There is a connection between this work and hitting times of stationary processes, or synonymously, the subject of escape rates in dynamical systems with holes. For a given stationary process, one is often concerned with the time taken to reach a particular state. A well-known instance of this is Penney's game \cite{Penney}, concerning the hitting time of certain sequences upon tossing of a fair coin, and the related `ABRACADABRA problem'. In the language of discrete dynamical systems, one has a phase space $X$ and a self-map $T:X \to X$. For each $x \in X$, one forms an orbit of iterates $x, T(x), T(T(x)), T(T(T(x))), \ldots $.  The system has a hole $H \subset X$ if $T$ is defined only on $X\setminus H$, so if the orbit ever falls into $H$ it cannot continue any further.  This hole can be thought of as a state of a certain stationary process where the starting point $x$ is chosen according to an invariant probability measure. In the context of our paper, one can start with an ambient shift space $X$. Forbidding a word $w$ from $X$ is similar to creating a hole for the cylinder set corresponding to $w$. In this case, the rate of ``escape'' into the hole is given exactly by the difference in entropy between the ambient shift and the subshift with the hole removed. One can also consider a system with multiple holes.

This line of research traces back to \cite{MR0534126} for expansive maps. Subsequently this has been studied extensively by several authors in a large body of work: \cite{MR2783972} (for systems like rotations, expanding maps on the interval and shift spaces), \cite{MR2535206} (for expanding maps and any absolutely continuous invariant measure),  \cite{MR2995652} (for Gibbs measures), \cite{MR4615466} (for one dimensional shifts of finite type) and \cite{agarwal2024escape} (for Markov chains). One key difference between these works and ours is the following: In a lot of these papers, they study the `hitting time/escape rate' of the hole given an ambient system. In our case, we are more interested in studying the dynamics on the space where all orbits touching the hole $H$ are removed, that is, instead of studying the dynamical system $(X, T)$ with an emphasis on the visits to the hole $H$ we look at the system $(X\setminus(\cup_{i\in \Z}T^{i}(H)), T)$. For more about these directions of research we refer the reader to the given citations and the references within.

%

%
%

%
%
\section*{Acknowledgements}
The first author would like to thank Aditya Thorat for many discussions and brainstorming ideas around the problems. He would like to thank Sharvari N. Tikekar and Haritha Cheriyath for first introducing these problems to him. He was supported by INSA Associate fellowship and the SERB Startup research grant. The third author was partially supported by ERC Synergy under Grant No. 810115 - DYNASNET and ERC Consolidator Grant 772466 “NOISE". We also thank Sophie MacDonald for helpful discussions and acknowledge support from NSERC.
\section{Earlier results in one dimension }\label{section:earlier_1d}

Some of the notations and definitions used in this section, such as full shift,
shift of finite
type (SFT), the SFT $X_{\{w\}}$, and entropy $h(X_{\{w\}})$ were given  in the introduction.

A finite directed graph $G$ is called {\em irreducible} if for every
ordered pair of vertices $I,J$ of $G$ there is a directed path from $I$ to $J$ (and otherwise
called {\em reducible}). An SFT $X$ is called {\em irreducible} if for any two allowed words $u, v$, there is an allowed word $w$ such that $uwv$, the concatenation of these three words, is also allowed in $X$. Otherwise, $X$ is a {\em reducible} SFT.

A {\em presentation} of $X_{\{w\}}$ is a finite directed graph $G$, together  with a labeling
of its edges,
such that $X_{\{w\}}$ is the set of all bi-infinite label sequences of bi-infinite walks in $G$. It is well-known that $X_{\{w\}}$ is irreducible if and only if it has an irreducible presentation (see, for example \cite[Proposition 3.3.11]{LM21}). 
A presentation
is called a {\em conjugacy presentation} if the mapping from bi-infinite edge sequences
to bi-infinite label
sequences in $X_{\{w\}}$ is one-to-one (and hence a bijection). This mapping is a sliding block
code with $m=a=0$ and thus is a conjugacy between $X_{\{w\}}$ and the edge shift $X_G$, hence the
name conjugacy presentation.
An example of a conjugacy presentation for any $X_{\{w\}}$ is
the follower set graph described at the beginning of Section~\ref{section:follower}.

For an irreducible graph $G$, the adjacency matrix $A_G$ has a unique maximum (in modulus)
eigenvalue  $\lambda_G$ such that $\lambda_G>0$, called the {\em Perron eigenvalue}.  It has
a unique right (Perron) eigenvector and a unique left (Perron) eigenvector, up to scale, and they are
both strictly positive vectors. Let $\chi_G(t) = t^{-d}p_G(t)$ where $p_G(t)$ is the characteristic
polynomial of $A_G$ and $d$ is the unique integer such that $\chi_G(t)$ has a nonzero constant term.
One can show that all graphs $G$ which have a labeling that is a conjugacy presentation of
$X_{\{w\}}$ have the same $\lambda_G$ and the same $\chi_G(t)$ which we shall respectively denote by $\lambda_w$
and $\chi_w(t)$, and call the Perron eigenvalue and characteristic polynomial of
$X_{\{w\}}$. ($\lambda_w$ is the largest root of $\chi_w(t)$.)
For more information on Perron-Frobenius theory, see \cite[Chapter 4]{LM21}.

The {\em entropy} of a shift space $X$ is defined as
$$
h(X) = \lim_{n\to\infty} (1/n)\log |\euscript{B}_n(X)|
$$
where $\euscript{B}_n(X)$ is the set of words of length $n$ that occur in elements of $X$.
In the case of a single forbidden word $w$, $\euscript{B}_n(X)$ was denoted by $\euscript{B}_n(w)$
in the introduction, and
$$h(X_{\{w\}}) = \lim_{n\to\infty} (1/n)\log |\euscript{B}_n(w)|.$$
It is well known that $h(X_{\{w\}}) = \log \lambda_w$.

We now define the correlation polynomial, which we mentioned in the introduction.
For any positive integer $q$ and any word $w=w_1\cdots w_k$ of length $k$ over the alphabet $[q]$, define the {\em self-overlap set} of $w$ by
$$
\ov(w,w):= \{i\in \{1, 2, \ldots, k\}: w_1 w_2 \cdots w_i = w_{k-i+1}w_{k-i+2} \cdots w_k\}
$$
and the {\em correlation polynomial} of $w$ by
$$
\phi_w(t)=\sum_{i\in \ov(w,w)} t^{i-1}.
$$
For instance $\phi_{100010}(t) =  t^5 + t $.
It can be derived from \cite[Theorem 1]{Lin89} that
\begin{align} \label{Doug_formula}
\chi_w(t)=(t-q)\phi_w(t)+1.
\end{align}

Let $p_n(X_{\{w\}})$ be the number of periodic points of $X_{\{w\}}$ with period $n$.
The {\em zeta function} of $X_{\{w\}}$  is defined as
$\zeta_{w}(t)=\exp\left( \sum_{n=1}^\infty \frac{p_n(X_{\{w\}})}{n} t^n\right)$.
It is known that for any conjugacy presentation of  $X_{\{w\}}$,
$\zeta_{w}(t)=\frac{1}{\det(Id-tA_{G})}$ where $Id$ is the identity matrix
having the same size as $A_G$ (see, for example, \cite[Theorem 6.4.6]{LM21}). Thus, it can be
readily seen that $\chi_w(t)$ determines $\zeta_{w}(t)$.
It follows that the information contained in $\zeta_w(t)$ is exactly
the set of non-zero eigenvalues of $A_G$.

We use $\tau_w$ to denote the {\em hitting time} random variable of $w$ by an i.i.d.\ 
sequence of Uniform($[q]$) random variables $(Y_i)_{i\in \mathbb{N}}$, i.e.,
$$
\tau_w=\min\{t\in \mathbb{N}: Y_{t-k+1} Y_{t-k+2} \cdots Y_t =w\}.
$$

Earlier results on the relations between quantities introduced above are summarized in the 
following two propositions. We made an effort to present these results in an organized way, and to give straightforward proofs of the various implications.


\begin{pr} \label{equality}
	Let $u, v$ be two strings of length $k$ over the $q$-ary alphabet ($q\geq 2$). Then the following are equivalent:
	\begin{enumerate}
		\item[(1)] $\mathrm{overlap}(u,u)=\mathrm{overlap}(v,v)$;
		
		\item[(2)] $\phi_u(t)=\phi_{v}(t)$ for all $t$;
		
		\item[(3)] $\phi_u(q)=\phi_{v}(q)$;
		
		\item[(4)] $E\tau_u=E \tau_{v}$;
		
		\item[(5)] $\zeta_{u}(t)=\zeta_{v}(t)$ for all $t$;
		
		\item[(6)] $\vert\euscript{B}_n({u})\vert=\vert\euscript{B}_n({v})\vert$ for all $n$;
		
		\item[(7)] $\vert\euscript{B}_n({u})\vert=\vert\euscript{B}_n({v})\vert$ for all sufficiently large $n$;
		
		\item[(8)] $h(X_{\{u\}})=h(X_{\{v\}})$.
	\end{enumerate}
\end{pr}
\noindent{\bf Proof:}
{{\em (1)} $\iff$ {\em (2)} $\iff$ {\em (3)}}: {\em (1)} $\Rightarrow$ {\em (2)} $\Rightarrow$ {\em (3)} is obvious. To see {\em (3)} $\Rightarrow$ {\em (1)}, note that for any $w$, the coefficients of $\phi_w(t)$, when viewed as a string, is the $q$-ary expansion of $\phi_w(q)$. Thus, $\phi_u(q)=\phi_v(q)$ implies $\phi_u(t)$ and $\phi_v(t)$ are the same polynomial, which in turn gives {\em (1)} by the definition of the correlation polynomial.

\medskip

\noindent{\bf {\em (3)} $\iff$ {\em (4)}}: This follows from the fact that 
$E\tau_w = q\phi_w(q)$ for any string $w$ (See \cite[pg. 189]{GO81},\cite{Sov66},\cite{Li83}).
\medskip

\noindent {\bf  {\em (1)} $\Rightarrow$ {\em (6)}}:
 This follows from \cite[pg. 185]{GO81}, where they showed that for any word $w$, the generating function of $\euscript{B}_n({w})$ is determined by $\phi_w(t)$. More precisely,
		$$
		G_w(t) = \frac{t \phi_w(t)}{1-(t-q) \phi_w(q)},
		$$
		where $G_w(t)$ is the generating function defined at the beginning of the introduction.
	
\noindent {\bf {\em (6)} $\Rightarrow$ {\em (7)} $\Rightarrow$ {\em (8)}}: Clear from the definition of entropy.

\noindent {\bf {\em (8)} $\Rightarrow$ {\em (1)}}:  The ideas in this proof are largely in
Eriksson
~\cite{Eri97}. It suffices to show {\em (8)} $\Rightarrow$ {\em (3)}. We prove the contrapositive.
	Assume WLOG that $\phi_u(q)>\phi_v(q)$. We now show $\lambda_u>\lambda_v$. To this end, we first write
	$$
	\phi_u(t)=\sum_{i=0}^{k-1} a_i t^i, \quad \phi_v(t)=\sum_{i=0}^{k-1} b_i t^i
	$$
	where $a_i, b_i \in \{0,1\}$. Define $r:= \max\{i: a_i=1, b_i=0\}$. Note that $r\leq k-2$ and $\phi_u(t)-\phi_v(t) > t^r-\sum_{i=0}^{r-1} t^i$ for all $t> 0$.
	Let $\lambda^*$ be the largest {positive} root of $t^r-\sum_{i=0}^{r-1} t^i=0$.

	We claim that $\lambda_v>\lambda^*$. To see this, let $w$ be a word of length $k$ with trivial self-overlap, i.e., $\phi_w(t)=t^{k-1}$. Since $\phi_v(t)\geq \phi_w(t)$ for all $t> 0$, we first deduce from equation (\ref{Doug_formula}) that $\chi_v(\lambda_w)=\chi_v(\lambda_w)-\chi_w(\lambda_w)\le 0$ and therefore $\lambda_v\ge\lambda_w$. It remains to show $\lambda_w>\lambda^*$. To this end, first observe that when $q=2$, the characteristic polynomial of $X_{\{w\}}$ is $\chi_w(t)=(t-2)t^{k-1}+1= (t-1)(t^{k-1}-\sum_{i=0}^{k-2} t^i)$. Thus, when $q=2$, $\lambda_w$ is the largest root of $t^{k-1}-\sum_{i=0}^{k-2} t^i=0$ and therefore $\lambda_w>\lambda^*$ (because the largest root of $t^s-\sum_{i=0}^{s-1}t^i=0$ is an increasing function of $s$). Finally, noting from equation (\ref{Doug_formula}) that the value of $\lambda_w$ increases as $q$ increases,
we have $\lambda_w>\lambda^*$ for all $q\geq 2$ and the claim is proven.

	Now,  since $\lambda_v<q$, we infer from equation (\ref{Doug_formula}) that
	$$
	\chi_u(\lambda_v)-\chi_v(\lambda_v)=(\lambda_v-q)[\phi_u(\lambda_v)-\phi_v(\lambda_v)]<(\lambda_v-q) [(\lambda_v)^r-\sum_{i=0}^{r-1} (\lambda_v)^i].
	$$
Since $q>\lambda_v>\lambda^*$, we have $\chi_u(\lambda_v)=\chi_u(\lambda_v)-\chi_v(\lambda_v)<0$ (recalling that $\chi_v(\lambda_v)=0$) and therefore $\lambda_u>\lambda_v$.

\noindent{\bf {\em (5)} $\iff$ (8)}: {\em (5)} $\Rightarrow$ {\em (8)}; this
follows from the fact
that the entropy is determined by the largest positive eigenvalue
 of $A_G$ above and the zeta function is determined by all of the non-zero eigenvalues of $A_G$.
To show {\em (8)} $\Rightarrow$ {\em (5)}, it suffices to show {\em (2)} $\Rightarrow$ {\em (5)}. But this follows from equation (\ref{Doug_formula}) and the fact that $\chi_w(t)$ determines $\zeta_{w}(t)$. \qed



\medskip

Some of the equivalent equalities in Proposition \ref{equality} can be strengthened to inequalities:

\begin{pr} \label{inequality_1d}
	Let $u, v, q$ be defined as in Proposition \ref{equality}. Then, the following are equivalent:
	\begin{enumerate}
		\item[(A)] $\phi_u(q) > \phi_v(q)$;
		\item[(B)] $E \tau_{u}> E\tau_v$;
		\item[(C)] $\vert\euscript{B}_n({u}) \vert\geq \vert\euscript{B}_n({v})\vert$ for all $n$ and $\vert\euscript{B}_n({u})\vert>\vert\euscript{B}_n({v})\vert$ for all sufficiently large $n$;
		\item[(D)] $h(X_{\{u\}})>h(X_{\{v\}})$.
	\end{enumerate}
\end{pr}

\noindent{\bf Proof:}
{{\em (A)} $\iff$ {\em (B)}}: This follows again from the fact that 
$E\tau_w = q\phi_w(q)$ for any string $w$ \cite[pg. 189]{GO81}.
\medskip

\noindent{\em (A) $\Rightarrow$ (D)}: By the result in Section 7 of \cite{GO81}, $(A)$ implies $\vert\euscript{B}_n({u}) \vert\geq \vert\euscript{B}_n({v})\vert$ for all $n$. Thus, $h(X_{\{u\}})\geq h(X_{\{v\}})$. By Proposition \ref{equality} (the equivalence between item $(3)$ and item $(8)$), we must have $h(X_{\{u\}})>h(X_{\{v\}})$.
\medskip

\noindent{\em (D) $\Rightarrow$ (A)}: Suppose to the contrary that $\phi_u(q)\leq \phi_v(q)$, i.e., either $\phi_u(q) < \phi_v(q)$ or $\phi_u(q)=\phi_v(q)$. In the former case, we have $h(X_{\{u\}})<h(X_{\{v\}})$ (by $(A)\Rightarrow (D)$), a contradiction; in the latter case, we have $h(X_{\{u\}})=h(X_{\{v\}})$ (by the equivalence between item $(3)$ and item $(8)$ in Proposition \ref{equality}), which is again a contradiction.
\medskip

\noindent{\em (A) $\Rightarrow$ (C)}: The first part is the result in Section 7 of \cite{GO81}. To see the second part, note that $(A)$ implies $(D)$, which immediately gives $\vert\euscript{B}_n({u})\vert>\vert\euscript{B}_n({v})\vert$ for all sufficiently large $n$.

\medskip

\noindent{\em (C) $\Rightarrow$ (A)}: Suppose to the contrary that $\phi_u(q)\leq \phi_v(q)$, i.e., either $\phi_u(q) < \phi_v(q)$ or $\phi_u(q)=\phi_v(q)$. In the former case, we have $\vert\euscript{B}_n({u})\vert<\vert\euscript{B}_n({v})\vert$ for all sufficiently large $n$ (by $(A) \Rightarrow (C)$), a contradiction; in the latter case, we have $\vert\euscript{B}_n({u})\vert=\vert\euscript{B}_n({v})\vert$ (by the equivalence between item $(3)$ and item $(6)$ in Proposition \ref{equality}), which is again a contradiction. \qed



\section{Graphical representation for 1D SFTs with one forbidden word}
\label{section:follower}

%

In this section we study a specific graph presentation for one dimensional SFTs $X_{\{w\}}$.
The construction allows direct comparison of different such SFTs, using ideas from both
probability, namely a coupling argument in Theorem \ref{pr:wine}, and from linear algebra,
namely the Perron Frobenius theorem in Theorem \ref{pr:cheese}. Together, these two theorems
enrich the equivalence $(B) \iff (D)$ from Proposition \ref{inequality_1d} under the
natural assumption \ref{assumeD}.


\subsection{Ambient shift is the full shift}
 \begin{figure}
 \centering
 \includegraphics[width = .66\textwidth]{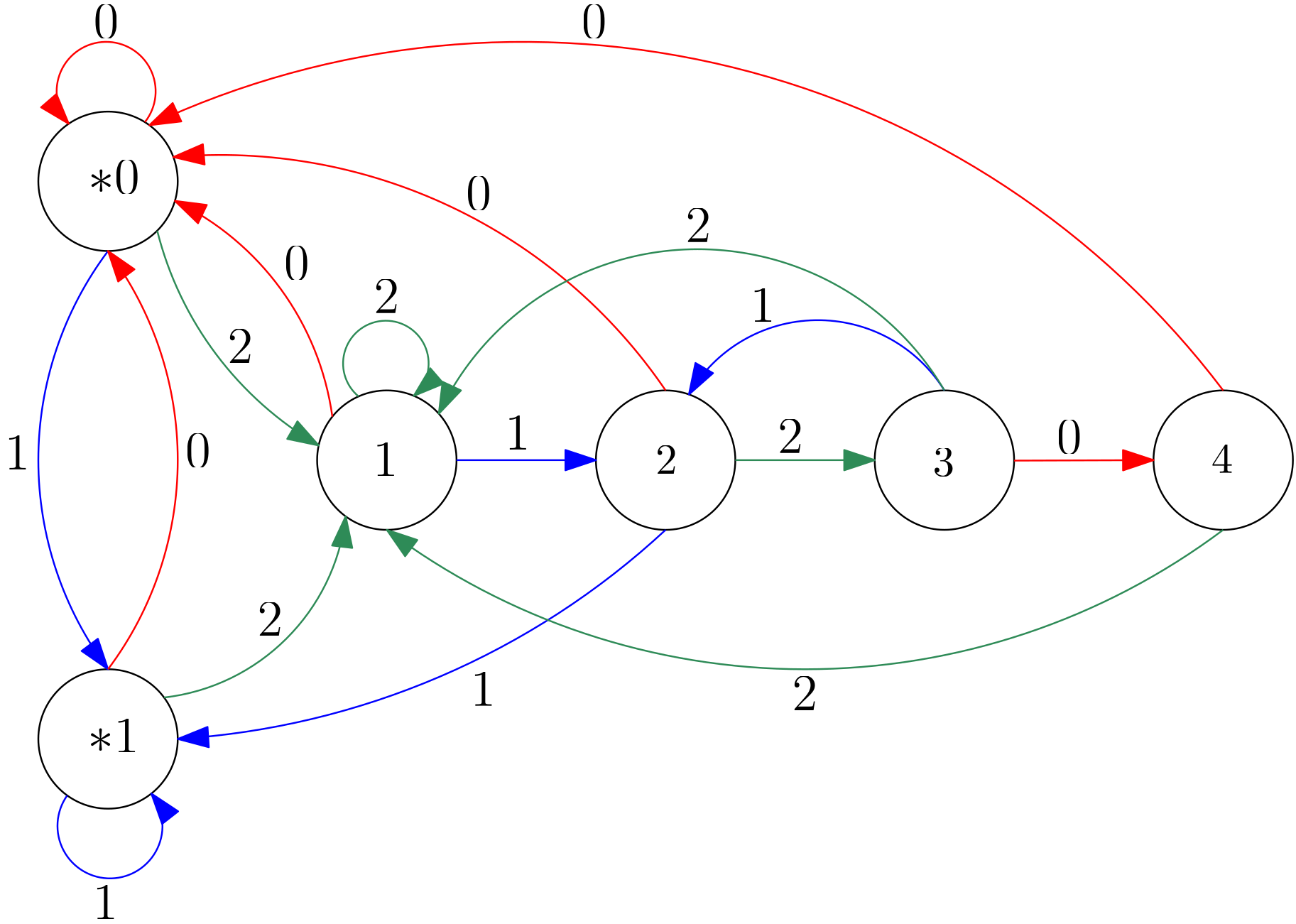}
 \caption{\label{figure: follower}The graph $L_{21201;3}$ for the word $w = 21201$ over the alphabet $[3] = \{0,1,2\}$. Edges labeled $0, 1, 2$ are colored in red, blue, and green, respectively.}
 \end{figure}

Fix $q \in \{2, 3, \ldots\}$, $k \geq 2$, and a word $w = w_1 w_2 \cdots w_k\in [q]^k$. First we introduce a bit of notation:

\begin{definition} \label{def:di} For any $a \in [q]$ and $1 \leq i \leq k-1$, define

\begin{equation} \label{didef} d_i(a) := \max\{j: w_1 w_2 \cdots w_{j} = w_{i-j+2} w_{i-j+3} \cdots w_i a\}, \end{equation}

with $d_i(a) = *a$ if the set of such $j$ is empty. \end{definition}

The $d_i$ can be viewed as a kind of self-overlap data: namely, $d_i(a)$ is the length of the longest prefix
of $w$ that is a suffix of the word $w_1 w_2 \cdots w_i a$.


\begin{definition} \label{def:fsg} Let $L_{w;q}$ be the following directed, edge-labeled graph on the vertex set $\{*a: a \in [q] \setminus \{w_1\}\} \cup \{1, 2, \ldots, k-1\}$. The edge set consists of:

\begin{itemize} \item For each $a, b \in [q] \setminus \{w_1\}$, an edge $*a \to *b$ labeled $b$;
\item For each $a \in [q] \setminus \{w_1\},$ an edge $*a \to 1$ labeled $w_1$;
\item For each $a \in [q]$ and $1 \leq i \leq k-1$, except for $a = w_k$ and $i = k-1$, an edge $i \to d_i(a)$ labeled $a$.
\end{itemize}

\end{definition}

For brevity we often suppress the dependence on $q$ and simply write $L_w$. An example with $q = 3$ is given in Figure \ref{figure: follower}.

 The idea behind this construction is the following: imagine a
 one-sided infinite sequence $x = (\ldots, x_{-2}, x_{-1}, x_0) \in [q]^\mathbb{N}$ that we
 modify one step at a time by adding a new letter on the right. The vertices of $L_w$ represent
 the length of the maximal suffix of $x$ that is a prefix of $w$, and the edges of $L_w$ represent
 appending that edge's label, say $a$, to $x$ to obtain $x' = (\ldots, x_{-1}, x_0, a)$.

The graph $L_w$ is, up to a small modification (merging the states $*a$ together), the
 {\em{follower set graph}} $G_w$ of the shift $X_{\{w\}}$, a well-known object in symbolic
 dynamics~\cite[Section 3.2]{LM21}. The follower set graph, described in the proof of
 Fact~\ref{follower} below,
   is a standard way to present any SFT (or even any sofic shift, a generalization of
SFTs).
We stick with the graph $L_w$, rather than $G_w$,
because a key lemma, Fact \ref{lem:graphtoword}, does not hold for the latter.
%

\begin{fact}
\label{follower} 
$L_w$ is a conjugacy presentation of $X_{\{w\}}$. That is, the mapping from bi-infinite edge sequences
to bi-infinite label sequences is a bijection onto $X_{\{w\}}$.
\end{fact}
\begin{proof}


Say that a word $u$ is {\em allowed} if it does not contain $w$.
For an allowed word $u$, define the {\em follower set}

\begin{equation}
F(u) = \{\text{words } v: uv \mbox{ is allowed}\}.
\end{equation}
We claim that there are only finitely many follower sets, in particular that
$F(u) = F(u^*)$ where $u^*$ is
the longest suffix of $u$ that is a proper prefix of $w$. Clearly $F(u) \subseteq F(u^*)$. For the reverse inclusion, let $v \in F(u^*)$.  If $v\not\in F(u)$,   then $w$ is contained in $uv$ but not in $u^*v$.  This means that $u$ has a longer (than $u^*$) suffix that is a proper prefix of $w$, a contradiction.

The {\em follower set graph} $G_w$ is a labeled directed graph whose vertices are $\{F(u)\}$ over all allowed $u$, and for each $a \in [q]$ s.t. $ua$ is allowed,
there is an edge, labeled $a$, from $F(u)$ to $F(ua)$. Note that by the above, $G_w$ is a finite labeled directed graph.

We claim that $G_w$ is a conjugacy presentation of $X_{\{w\}}$.
It is clear that the set of label sequences of $G_w$ is contained in $X_{\{w\}}$.

Observe that any edge sequence in $G_w$ with label sequence $x = \ldots ,a_n, a_{n+1}, \ldots$ is of the form $ = \ldots F(u^n) F(u^{n+1})  \ldots$ where each
$F(u^{n+1}) = F(u^na_n)$.  But then for fixed $n$, $F(u^{n+k}) = F(u^na_n \ldots a_{n+k}) = F(u^*)$ where $u^*$ is the smallest prefix of $w$ that is a suffix of $u^na_n \ldots a_{n+k}$; but
clearly $u^*$ is contained in  $a_n \ldots a_{n+k}$ and thus  $F(u^{n+k}) = F(u^*)$ can be uniquely recovered from $x$. Thus, every element of $X_{\{w\}}$ is the label sequence of a unique edge sequence in $G_w$, and so $G_w$ is indeed a conjugacy presentation of $X_{\{w\}}$.

Finally, we see that $L_w$ is an {\em in-splitting} of $G_w$ (for the definition, see ~\cite[Definition 2.4.7]{LM21}): identify each vertex $F(w_1 \ldots w_i)$ in $G_w$ with the vertex $i$ in $L_w$; in $G_w$, the vertex $F(\emptyset)$ has $q-1$ self-loops, each labeled by a unique element of $[q]\setminus \{w_1\}$, and, to form $L_w$ from $G_w$, that vertex is split  into $q-1$ vertices, denoted
$*a$, $a \in [q]\setminus \{w_1\}$.  Since in-splittings define conjugacies (\cite[Theorem 2.4.10]{LM21}), $L_w$ is also a conjugacy presentation of $X_{\{w\}}$.
%
%
%
%
%
%
\end{proof}




We now note some nice properties of $L_w$ which will be useful throughout this section. By convention, we write $*a < i$ for any $a \in [q]$ and $1 \leq i \leq k-1$.

\begin{pr} \label{prop:donkey_kong} Fix $w \in [q]^k$ for $k \geq 2$ and let $L_w$ and $d_i$
be as in Definitions \ref{def:di} and \ref{def:fsg}.
	
	\begin{enumerate}[i.]
		\item All vertices $i \leq k-2$ and $*a$ for $a \in [q] \setminus \{w_1\}$ have out-degree $q$, and state $k-1$ has out-degree $q-1$.
		
		\item Incoming edges to a vertex $i \geq 1$ have label $w_i$, and incoming edges to vertex $*a$ have label $a$.
	
		\item For all $i \in \{1, \ldots, k-2\}$ and $a \neq w_{i+1}$, $d_i(a) \leq i$; and $d_i(w_{i+1}) = i+1$.
	
		\item If $d_i(a), d_{i'}(a') \geq 1$, then $i - d_i(a) = i' - d_{i'}(a')$ only if $i = i'$, $a = w_{i+1}$ and $a' = w_{i'+1}$.
		
		\item Any prefix $u = w_1 w_2 \cdots w_r$ of $w$ with $r \geq 2$ has graph $L_u$ identical to the induced subgraph of $L_w$ on the vertices
		\begin{equation} \{*a: a \in [q] \setminus \{w_1\}\} \cup \{1, 2, \ldots, r-1\} \end{equation}
		
	\end{enumerate}
	
\end{pr}


\begin{proof}
Parts $i$, $ii$, and $v$ are clear from the definition. For part $iii$, it follows from the definition of $d_i(a)$ that $d_{i}(a) \leq i+1$ for all $a$, and from the edge-labeling that $d_i(a) \neq d_i(a')$ if $a \neq a'$. For part $iv$, write $d_i = d_i(a)$ and $d_{i'} = d_{i'}(a')$ for short, and suppose for the sake of contradiction that $i - d_i = i' - d_{i'}$ for some $i \neq i'$ with $i, i' \geq 1$. Then also $d_i \neq d_{i'}$ by re-arranging, so assume WLOG $d_i > d_{i'}$. If $a = w_{i+1}$, then by part $iii$, $a' = w_{i'+1}$, and similarly with the roles of $a$ and $a'$ reversed, so assume $a \neq w_{i+1}$ and $a' \neq w_{i'+1}$. Unraveling the definition of $d_i$ gives
\begin{equation} w_{i-d_i+j+1} = w_j \text{ for } j = 1, 2, \ldots, d_{i}-1, \text{ and } w_{i+1} \neq w_{d_i}. \end{equation}
Since $d_i> d_{i'} \geq 1$ by assumption, setting $j = d_{i'}$ gives
\begin{equation} w_{d_{i'}} = w_{i - d_{i} +d_{i'} + 1} = w_{i' - d_{i'} + d_{i'} + 1} = w_{i' + 1}, \end{equation}
contradicting the definition of $d_{i'}$.

\end{proof}


We remark that the properties in Proposition \ref{prop:donkey_kong} are not sufficient to
characterize the set of graphs that arise as $L_w$ for some word $w$ and integer alphabet
size $q \geq 2$, and we are not aware of such a set of sufficient conditions. We now give
a few applications of Proposition \ref{prop:donkey_kong}, starting with a characterization
of which of the $L_w$ graphs are irreducible. Recall that a graph is {\em irreducible} if
for any vertices $x, y$, there is a (directed) path from $x$ to $y$ (and otherwise it is {\em reducible}).

\begin{pr} \label{lem:irr} The graph $L_w$ is reducible if and only if $q = 2$ and $
w$ is one of the words $10^{k-1}, 1^{k-1}0, 01^{k-1}, 0^{k-1}1$ for some $k \geq 2$.
\end{pr}

\begin{proof} Since the graph $L_w$ always has the paths $*a \to 1 \to \cdots \to k-1$ for $a \in [q] \setminus \{w_1\}$ and contains an induced complete graph on the vertices labeled $*a$, it suffices to show that for each $i \geq 1$, $d_i(a) < i$ for some letter $a \in [q]$. We start with the easier case $q \geq 3$. Proposition \ref{prop:donkey_kong} parts $i$ and $iii$ together imply that there is exactly one edge $i \to i+1$ and at most one edge $i \to i$, and thus at least one edge $i \to d_i(a) < i$.

So assume $q = 2$. Write $\overline{a} = 1-a$ and $d_i = d_i(\overline{w_{i+1}})$ for the binary alphabet.  Suppose that $L_w$ is reducible, i.e. that there is no path $k-1 \to *a$ for any $a$. Then there exists $1 \leq i \leq k-1$ such that $d_j \geq i$ for all $j = i, i+1, \ldots, k-1$: let $I$ be the largest such $i$. If $I = k-1$, then $d_{k-1} = k-1$, and unraveling the definition of $d_{k-1}$ gives $w_1 = w_2 = \cdots = w_{k-1} = \overline{w_k}$, i.e. $w = 1^{k-1}0$ or $0^{k-1} 1$.

So suppose $I < k-1$. By part $iv$ of Proposition \ref{prop:donkey_kong}, since $I \geq 1$, the values $j - d_j$ must all be distinct for $j \in \{I, I+1, \ldots, k-1\}$. In particular, by part $iii$ of Proposition \ref{prop:donkey_kong}, $d_{I} = I$, and by induction over $j$, $d_j = I$ for all $j \in \{I, I+1, \ldots, k-1\}$. Writing out the definition of $d_j$, we obtain that for $j \in \{I, I+1, \ldots, k-1\}$,

\begin{equation} \label{eq:lucky} w_m = w_{j+1-I+m} \text{ for } m \in \{1, \ldots, I-1\}, \text{ and also } w_{I} = \overline{w_{j+1}}. \end{equation}

If $I \geq 2$, then taking $m = I-1$ with $j = I$ and separately with $I+1$ in Equation \ref{eq:lucky} immediately gives the contradiction $w_{I-1} = w_{I} = w_{I+1}$ and $w_{I}= \overline{w_{I+1}}$. Thus $I = 1$ and $w_1 = \overline{w_2} = \overline{w_3} = \cdots = \overline{w_{k}}$, i.e. $w = 10^{k-1}$ or $01^{k-1}$, as desired.

\end{proof}

\begin{rem}\label{Remark: Reducible subshift}
Using Proposition \ref{lem:irr}, one can easily show that the corresponding shift space
$X_{\{w\}}$ is reducible in the same cases, and irreducible otherwise. \end{rem}

Next we show that the graph $L_w$ is essentially enough information to determine the word $w$. In what follows, we say a map $\varphi$ between graphs $L_w, L_{w'}$ is a {\em{graph isomorphism}} if it is an isomorphism of the underlying unlabeled graphs, regardless of how $\varphi$ transforms the edge/vertex labels.

\begin{fact} \label{lem:graphtoword}
Let $q \geq 2$. If $w, w' \in [q]^k$ are two words such that $L_w, L_{w'}$ are isomorphic by an isomorphism $\varphi$ that preserves the vertex names, then $w = w'$. ($\varphi$ `preserving the vertex names' means: viewed as a map between the vertex sets, $\varphi(i) = i$ for all $i$ and $\varphi(*a) = *a$ for all $*a$ vertices.)
\end{fact}

\begin{proof} It suffices to show that given the graph $L_w$ with all edge labels erased, one can recover the
word $w$. We give a simple recursive algorithm to recover $w$ along with the edge labels of
$L_w$. For each vertex $*a$ that appears in $L_w$, assign all the incoming edges
to it the label $a$. By the construction of $L_w$, $w_1$ must be the unique digit in $[q]$ not
appearing among the $*a$ vertices. Assume by induction we have determined the values
$w_1, w_2, \cdots, w_j$ for some $1 \leq j \leq k-1$ and the labels of all incoming edges to vertices $1, 2, \ldots, j$. Then by the construction of $L_w$, $w_{j+1}$ is the unique digit in $[q]$ not appearing in the set of labels of the outgoing edges from vertex $j$, which by Proposition \ref{prop:donkey_kong} parts $ii$ and $iii$ and induction are among the known values. \end{proof}

One can easily strengthen Fact \ref{lem:graphtoword} to the following:

\begin{thm}\label{theorem: isomorphism of graphs}
Let $q\geq 2$. If $w, w' \in [q]^k$ are such that $L_w$ and $L_{w'}$ are isomorphic, then there is a permutation $\pi$ of the alphabet $[q]$ such that $\pi(w_i) = w'_i$ for $i = 1, 2, \ldots, k$.
\end{thm}
Theorem \ref{theorem: isomorphism of graphs} implies:
\begin{co}\label{recover_word} Let $q \geq 2$. For any words $w, w' \in [q]^k$, $L_w$ and $L_{w'}$ are isomorphic if and only if there is a permutation $\pi$ of the alphabet $[q]$ such that $\pi(w_i) = w'_i$ for $i = 1, 2, \ldots, k$.  \end{co}

	For example, in the binary alphabet case $q = 2$, the word $w$ can be recovered from the vertex labels on $L_w$ from the following simple recursive formula. In this case, we can write $*a$ for the unique $*$ vertex in $L_w$, so that $w_1 = \overline{a}$ (recalling the notation from the proof of Proposition \ref{lem:irr}). Then for each $i\in \{1, \ldots, k\}$, one easily finds
	
	$$w_{i}=
	\begin{cases}
		w_1, & \text{ if } d_{i-1}=*a \\
		\overline{w_{d_{i-1}}}, & \text{ otherwise}
	\end{cases}
	$$

\begin{proof}[Proof of Theorem \ref{theorem: isomorphism of graphs}]
We will first show that the isomorphism of the graphs must preserve certain labels. Indeed if $\varphi$ is the isomorphism from $L_w$ to $L_{w'}$ then $\varphi(k-1)=k-1$. This is true since $k-1$ is the only vertex with degree $q-1$ in both the graphs. We will proceed by induction to show that $\varphi (t)=t$ for all $1\leq t \leq k-1$. Suppose we know this for $t\geq j$ for some $j$. Now the only incoming edge to the vertex labeled $j$ which is not amongst $\{j, j+1, j+2, \ldots, k-1\}$ must originate at $j-1$. Thus $\varphi (j-1)=j-1$. This completes the induction step and implies that $\varphi(t)=t$ for all $1\leq t\leq k-1$. Without loss of generality assume that $w_1=w'_1 = 0$. Then $\varphi$ is a bijection from vertices labeled $\{*s~:~1\leq s \leq q-1\}$ in $L_w$ to the same set of vertices in $L_{w'}$. This induces a bijection on $[q]$ which by Fact \ref{lem:graphtoword} takes $w$ to $w'$.
\end{proof}

By Fact \ref{follower}, an isomorphism of the graphs $L_w$ and $L_{w'}$ implies conjugacy of the corresponding shift spaces. As we will learn in Theorem \ref{conj_trivial_self_over} of Section \ref{section:conjugacy}, many pairs of
words $w, w'$, are not equivalent under a permutation of $[q]$ and yet give rise to conjugate shift spaces $X_{\{w\}}, X_{\{w'\}}$.

Also, we note that Corollary \ref{recover_word} implies that the correlation polynomial $\phi_w$ is uniquely determined by the graph $L_w$. Despite this, we do not know an explicit map between the family of possible correlation polynomials and the family of possible graphs $L_w$. See section \ref{sec:discussion} for further discussion.

We now turn to the main results of this section, which obtain the same conclusions as in ~\cite{GO81}, but under
a different assumption that is adapted to Perron Frobenius theory. We first make a short excursion into discrete probability. First recall the stochastic dominance order for random variables:



\begin{definition} For any random variables $A, B$, possibly defined on different probability spaces, write $A \prec_{\text{st}} B$ if $\mathbb{P}(A \geq t) \leq \mathbb{P}(B \geq t)$ for all $t \in \mathbb{R}$. Equivalently, $A \prec_{\text{st}} B$ if and only if there is a probability space $(\Omega, \mathbb{P})$ and a pair of random variables $(X, Y)$ on $(\Omega, \mathbb{P})$ such that $X$ and $Y$ have the distributions of $A$ and $B$ respectively, and $\mathbb{P}(X \leq Y) = 1$. \end{definition}

Recall the random variable $\tau_w$, the hitting time of word $w$ in an iid Uniform$([q])$ sequence.
It follows easily from the results of \cite[pg. 204]{GO81} that
$h(X_{\{w\}}) \geq h(X_{\{w'\}})$ if and only if
 $\tau_{w'} \prec_{\text{st}} \tau_{w}$.
In Theorems \ref{pr:wine} and \ref{pr:cheese} we prove the same result under the
following natural assumption on the edges $d_i$. Fix words $w, w' \in [q]^k$ and let
$L = L_w, L' = L_{w'}$  and $d_i, d_i'$ as in Definition \ref{def:di}.

\begin{assumption}\label{assumeD} For words $w, w' \in [q]^k$, we say condition $D(w, w')$ holds if for each $i \in \{1, 2, \ldots, k-1\}$, there exists a permutation $\pi_i:[q] \to [q]$ such that $d_i(a) \leq d_i'(\pi_i(a))$ for all $a \in [q]$. \end{assumption}

In the binary alphabet case $q = 2$, the assumption $D(w, w')$ is
equivalent to $d_i \leq d_i'$ for all $i$ (recalling the notation from the proof of Proposition \ref{lem:irr}). For example, $D(w,w')$ holds for the pair $w = 01010, w' = 01000$; also, $D(w, 0^k)$ always holds for $w \in [q]^k$.
Put simply, $D(w,w')$ says that the
edges of $L$ go further back towards the $*$ vertices than the edges of $L'$.

While this assumption is natural in the context of the graphs $L_w$, we do not know any explicit formulation of Assumption \ref{assumeD} in terms of autocorrelation polynomials. See Section \ref{sec:discussion} for further discussion.

We note that this type of ordering has already been considered in more detail in \cite{DS16},
under the name `skip-free' Markov chains, which in
particular contains a more general version of Theorem \ref{pr:wine}. We include a simplified proof here for the sake of completeness.

%
%

\begin{thm}\label{pr:wine} Let $w, w' \in [q]^k$ for $k \geq 3$ be any words with graphs
$L = L_w, L' = L_{w'}$. If $D(w, w')$ holds, then $\tau_{w'} \prec_{\text{st}} \tau_{w}$. Additionally, if for some $1 \leq i \leq k-1$ and $a \in [q]$ we have a strict inequality $d_i(a) < d_i'(\pi_i(a))$, then $\tau_w$ and $\tau_{w'}$ are not equal in distribution. \end{thm}

\begin{proof} Consider the simple random walk on the graph $G_w$ obtained from $L_w$ by identifying all vertices $\{*a: a \in [q] \setminus \{w_1\}\}$ to a single vertex labeled $*$, and adding an extra vertex $k$ is added at the end of $L_w$, i.e. there is an edge $k-1 \to k$ corresponding to hitting the word $w$. Let $X, X'$ denote simple random walk on the graphs $G_w, G_{w'}$ respectively, stopped at the first time $\tau$ (resp. $\tau'$) that $X$ (resp. $X'$) hits state $k$. That is, $X$ is a random sequence of vertices $X = (* = X_0, X_1, X_2, \ldots, X_{\tau_{k}} = k)$ such that: given $X_i$, $X_{i+1}$ a uniformly random neighbor of $X_i$ for each $1 \leq i < \tau_k$. (For completeness one can add outgoing edges from vertex $k$ corresponding to appending any digit to $w$, so $G_w$ is a conjugacy representation for the full shift, and simple random walk on $G_w$ is the corresponding maximum entropy Markov chain.) Note that $\tau, \tau'$ are distributed as $\tau_w$ and $\tau_{w'}$, respectively.

We couple the pairs $(X, \tau)$ and $(X', \tau')$ as follows. Generate $X$ and $X'$ along with a random function $\sigma$ recursively, initialized at $X(0) = X'(0) = \sigma(0) = *$. For $s \geq 1$, given $X(u), X'(u)$ and $\sigma(u)$ for $u \leq s$, let $X(s+1)$ have the distribution of an independent simple random walk step from $X(s)$ (in the graph $G$), and define $(X'(s), \sigma(s))$ by

\begin{equation*} (X'(s+1), \sigma(s+1)) = \begin{cases} (X'(s) + 1, \sigma(s) + 1) & \text{ if }X'(s) = X(s) \text{ and } X(s+1) = X(s) + 1 \\ (d'_{X'(s)}(\pi_{X'(s)}(a)), \sigma(s) + 1) & \text{ if } X'(s) = X(s) \text{ and } X(s+1) = d_{X(s)}(a) \\ (X'(s), \sigma(s)) & \text{ if } X'(s) > X(s) \end{cases} \end{equation*}

(Note that we interpret $d_i(a) = *$ if $d_i(a) = *b$ for some $b \in [q]$.) In words, when $X$ and $X'$ are at the same level, they move together (after translating by the permutation $\pi_X$), and when they are at different levels, $X$ moves independently while $X'$ remains frozen.  The (random) counter $\sigma(s)$ is the number of steps taken by $X'$ when $X$ has taken $s$ steps. The process ends at time $\tau = \min \{s: X(s) = k\} = \min \{s: X'(s) = k\}$.
By induction, using the assumption $D(w, w')$, and because the only simple path in both graphs from $i \to j$ with $i < j$ is the path $i \to i+1 \to \cdots \to j-1 \to j$, for all $s \geq 0$ we have $X(s) \leq X'(\sigma(s))$ and $\sigma(s) \leq s$.

It follows immediately from the construction that $X(\cdot)$ has the distribution of simple random walk on $G_w$ stopped on hitting vertex $k$, and $\tau$ has the distribution of $\tau_w$. Similarly, $X'(\gamma(\cdot))$ has the distribution of simple random walk on $G_{w'}$ stopped on hitting vertex $k$, where $\gamma(t) = \min \sigma^{-1}(t)$ are the jump times for $X'$, and thus $\sigma(\tau)$ has the distribution of $\tau_{w'}$. Combining this with the previous paragraph, $\tau \geq \tau'$ under this coupling. Additionally, $\tau > \sigma(\tau)$ if $\sigma(s) < s$ for some $s < \tau$, and the latter occurs with positive probability if $d_i(a) < d_i'(\pi_i(a))$ for some $a$ and $i$. \end{proof}


Theorem \ref{pr:wine} says that under condition $D(w, w')$, the stochastic ordering of the hitting times $\tau_w$ can be read off from the graphs $L_w$. It turns out that entropies can also be compared easily when condition $D(w,w')$ holds. Recall from the Perron-Frobenius Theorem that the adjacency matrix of an irreducible graph has a unique maximum (in modulus) real eigenvalue $\lambda > 0$ (the Perron eigenvalue), and the corresponding left and right eigenvectors $\ell$ and $r$ are both positive. Thus:
%

\begin{lem} \label{positive} For any irreducible graph $L_w$,
the corresponding left and right eigenvectors $\ell, r$ are 
strictly positive. For the reducible cases appearing in
Lemma~\ref{lem:irr} one can check by direct computation (omitted) that the Perron 
eigenvalue and corresponding left and right eigenvectors are strictly positive, with the 
exception that $\ell_{*a}$ and $r_{k-1}$ may be $0$. \end{lem}

We will use the following property of the right Perron eigenvector of $L_w$:

\begin{pr}\label{revec_decr} Fix $w \in [q]^k$. Let $A = A_w$ denote the adjacency matrix of
$L_w$, and write $r = r_w$ for the right eigenvector of $A$ corresponding to the
Perron-Frobenius eigenvalue $\lambda = \lambda_w$ of $A$. Then 
$0 < \lambda - q + 1 < 1$ and 
$r$ has strictly
decreasing entries. In particular, the entries of $r$ decrease at least exponentially:
for $i = 1, \ldots, k-2$,


\begin{equation}\label{rvec_dec_induction} r_{i+1} \leq (\lambda - q+1) r_i,\end{equation}

and for all $a \in [q] \setminus \{w_1\}$,

\begin{equation} \label{rvec_dec_base} r_{1} = (\lambda - q + 1) r_{*a}\end{equation}
\end{pr}

\begin{proof} First, for any $a \in [q] \setminus \{w_1\}$,

\begin{equation} \lambda r_{*a} = r_1 + \sum_{b \in [q] \setminus \{w_1\}} r_{*b}. \end{equation}

Since the right hand side does not depend on $a$, we have $r_{*a} = r_{*a'}$ for any 
$a, a' \in [q] \setminus \{w_1\}$. Thus

\begin{equation} \lambda r_{*a} = (q-1) r_{*a} + r_1, \end{equation}
which gives Equation \ref{rvec_dec_base}. Note that by positivity of 
$r$, $\lambda - q + 1 > 0$.  Combining this with $\lambda < q$ 
(since the shift space $X_{\{w\}}$ is a proper subshift of the full $q$-shift, which has
entropy $\log q$), 
we get $0 < \lambda - q + 1 < 1$.
We proceed by induction: assume that 
Equation \eqref{rvec_dec_induction} holds for $i = 1, \ldots, j-1$.
Then 

\begin{equation} \lambda r_j = \sum_{a \in [q]} r_{d_j(a)} \geq r_{j+1} + (q-1) r_j. \end{equation}

Here the equality uses that $r$ is a right eigenvector, and the inequality
follows from Proposition \ref{prop:donkey_kong} part $iii$, the induction hypothesis, and
the fact that $\lambda < q$. 
Re-arranging gives Equation \ref{rvec_dec_induction}. The same calculation, along with equation \ref{rvec_dec_base}, proves the base case $i = 1$.
\end{proof}


The proof actually shows something a bit stronger: namely that for $i = 1, \ldots, k-2$,

\begin{equation} \frac{r_{i+1}}{r_i} \leq \lambda - \sum_{a \in [q] \setminus \{w_{i+1}\} } (\lambda - q + 1)^{-i+d_i(a)}. \end{equation}

We can now describe the ordering of the  $h(X_{\{w\}})$  under the same assumption as in
Theorem \ref{pr:wine}.
%

\begin{thm} \label{pr:cheese} Fix $w, w' \in [q]^k$ for $k \geq 2$ and let
$r$ be as in Proposition \ref{revec_decr}. If
$D(w, w')$ holds, then $h(X_{\{w\}}) \geq h(X_{\{w'\}})$. If $d_i(a) < d_i'(\pi_i(a))$ for some $1 \leq i \leq k-2$ and $a \in [q]$, or if $L_w, L_{w'}$ are both irreducible and $d_{k-1}(a) < d_{k-1}(\pi_{k-1}(a))$ for some $a \in [q]$, then the entropy inequality is strict. \end{thm}

\begin{proof} Let $A, A'$ denote the adjacency matrices of $L, L'$, and write $r, r'$ for the
right eigenvectors of $A, A'$ corresponding to the top eigenvalues $\lambda, \lambda'$ (where all matrices and vectors are indexed by the vertex labels of the corresponding graphs). Recall that
  $h(X_{\{w\}}) = \log \lambda_w$ and
  $h(X_{\{w'\}}) = \log \lambda_{w'}$~\cite[Theorem 4.3.3]{LM21}. Combining Assumption $D(w,w,')$ and Proposition \ref{revec_decr}, we obtain that for all $i \in \{1, \ldots, k-2\}$

\begin{equation} \sum_{a \in [q]} r_{d'_i(a)} \leq \sum_{a \in [q]} r_{d_i(a)}, \end{equation}
and that

\begin{equation} \sum_{a \in [q] \setminus \{w_{k}\}} r_{d'_{k-1}(a)}
 \leq \sum_{a \in [q] \setminus \{w'_{k}\}} r_{d_{k-1}(a)}. \end{equation}
Thus for $i = 1, \ldots, k-2$,
\begin{equation} (A' r)_{i} = \sum_{a \in [q]} r_{d'_i(a)} \leq \sum_{a \in [q]} r_{d_i(a)} = (A r)_i = \lambda r_i, \end{equation}
and similarly, $(A' r)_{k-1} \leq \lambda r_{k-1}$. Also, we may assume without loss of generality (by applying a permutation of the alphabet) that $w_1 = w'_1$. Note that in this case the outgoing edges of the vertices $\{*a: a \in [q] \setminus \{w_1\}\}$ in $L$ and $L'$ have a bijective correspondence which preserves labels, so for all $a \in [q] \setminus \{w_1\}$ we have by an easy computation that $(A' r)_{*a} = \lambda r_{*a}$.

Combining these observations gives $A'r \leq \lambda r$. Take a left eigenvector $\ell'$ for the Perron-Frobenius eigenvalue $\lambda'$ in $A'$ and calculate:

\begin{equation} \lambda' \ell' \cdot r = \ell' A' r \leq \ell' (\lambda r) = \lambda \ell' \cdot r. \end{equation}

(Here $\cdot $ is the dot product.) The entries $\ell'_i$ and $r_i$ are strictly positive for $1 \leq i \leq k-2$, and by Lemma \ref{positive}, the first and last entries are non-negative. It follows that $\lambda' \leq \lambda$. By Lemma \ref{positive}, if some strict inequality $d_i(a) < d_i'(\pi_i(a))$ holds for $1 \leq i \leq k-2$, then $(A' r)_i < \lambda r_i$, the same argument yields $\lambda < \lambda'$; and the same conclusion follows if $L_w, L_{w'}$ are irreducible and $d_{k-1}(a) < d'_{k-1}(\pi_{k-1}(a))$ for some $a$.
\end{proof}

\subsection{The ambient shift is the golden mean shift} \label{fol_mem1}
In this section we generalize the $L_w$ graph construction and Theorem \ref{theorem: isomorphism of graphs} to the case where the underlying shift space is not the full $q$-shift, but the golden mean shift. Let us set things up first. Let $q\geq 2$ and $\mathcal F\subset [q]^2$. If $w\in \euscript B_n(X_{\mathcal F})$, we are broadly interested in learning the effect of forbidding a further word $w$. We will denote by $(X_{\mathcal F})_{\{w\}}$ the SFT with the forbidden set $\mathcal F\cup \{w\}$.
\begin{de} For $w \in \euscript{B}_n({\mathcal{F}})$ of length $k$, start with the graph $L_w$ exactly as in Definition \ref{def:di}. Define a new graph $L^{\mathcal{F}}_w$ by deleting edges as follows: for any path of length $2$ consisting of edges $e, f$ labeled $a_e, a_f$ with $a_e a_f \in \mathcal{F}$, delete edge $f$.
\end{de}

\begin{fact} The resulting graph $L^{\mathcal{F}}_w$ is
a conjugacy presentation of the shift space $(X_{\mathcal{F}})_{\{w\}}$. \end{fact}

This follows from Fact \ref{follower}.

Note that in the binary case $q = 2$, except for trivial shift spaces, we can take $\mathcal{F} = \{11\}$ without loss of generality. The graph $L^{\mathcal{F}}_w$ is identical to $L_w$ except that each vertex $i$ with $w_i = 1$ (and the vertex $*1$ if it is present in $L$) only has one outgoing edge $i \to i+1$. As a result, the word $w$ is determined in a simple way from the labeled graph $L^{\{11\}}_w$ (analogous to Fact \ref{lem:graphtoword}):

\begin{fact} \label{g2w2} Given the graph $L = L^{\{11\}}_w$ for an allowable word $w \in X_{\{11\}}$ of length $k$, for $i \in \{1, \ldots, k-2\}$,

\begin{equation} w_i = \begin{cases} 1, & \text{outdegree}_L(i) = 1 \\ 0, & \text{outdegree}_L(i) = 2, \end{cases} \end{equation}

and

\begin{equation} w_{k-1}w_k = \begin{cases} 10, & \text{outdegree}_L(k-1) = 0 \\ 0\overline{w_{d_{k-1}}}, & \text{outdegree}_L(k-1) = 1 \end{cases}  \end{equation} \end{fact}

With a bit more work, we also get the same conclusion as in Theorem \ref{theorem: isomorphism of graphs} for the golden mean shift.

\begin{thm} \label{isogen} For any $k \geq 3$ and $w, w' \in \euscript{B}_k({\{11\}})$,
$w = w'$ if and only if $L^{\{11\}}_w$ and $L^{\{11\}}_{w'}$ are isomorphic as unlabeled
graphs. \end{thm}

In order to prove this theorem, we need the following facts.
\begin{fact} \label{12cycle}
There is at most one $1$-cycle and at most one $2$-cycle in $L_w^{\{11\}}$.
\end{fact}
\begin{proof}
This follows immediately from the fact that $(X_{\{11\}})_{\{w\}}$ contains only one fixed point (i.e., $0^{\infty}$) and only one orbit of length $2$ (i.e., $\{(01)^\infty, (10)^\infty\}$), and the fact that $L_w^{\{11\}}$ is a conjugacy presentation of $(X_{\{11\}})_{\{w\}}$.
\end{proof}


\begin{fact} \label{1cycle}
If $L_{w}^{\{11\}}$ has a self-loop at some vertex $j$ where $j\geq 1$, then we must have $w_1 w_2\cdots w_j w_{j+1}= 0^j1$, i.e., the $(j+1)$-prefix of $w$ must be $0^j 1$.
\end{fact}

\begin{proof}
First note that the self-loop must be labeled $0$ (which means $\overline{w_{j+1}}=0$). Now, having a self-loop at vertex $j$ exactly means that the longest possible overlap between $w_1w_2w_3\cdots w_{j} \overline{w_{j+1}}$ and $w_1 w_2\cdots w_k$ is of length $j$, which in
particular implies $w_2w_3\cdots w_{j} \overline{w_{j+1}}=w_1 w_2\cdots w_j$. Thus, one immediately obtains $w_1=w_2=w_3=\cdots=w_j=\overline{w_{j+1}}=0$.
\end{proof}

\begin{fact} \label{2cycle}
If $L_{w}^{\{11\}}$ has a $2$-cycle at the two vertices $j, j+1$ where $j\geq 2$, then we must have
$$
w_1w_2\cdots w_j w_{j+1} w_{j+2}=
\begin{cases}
10101\cdots 01 00 \quad \mbox{if $j$ is odd} \\
010101\cdots 01 00 \quad \mbox{ if $j$ is even}.
\end{cases}
$$
\end{fact}
\begin{proof}
The edge from $j$ to $j+1$ must be labeled $0$. Thus, the $2$-cycle must be labeled $00$ or $01$. Since $L_w^{\{11\}}$ is a conjugacy presentation it follows that the $2$ cycle must be labeled $01$ (i.e., $w_{j+1} \overline{w_{j+2}}$=01).

Now, on the other hand, the existence of the $2$-cycle at $j,j+1$ implies that the longest possible overlap between $w_1w_2\cdots w_j w_{j+1} \overline{w_{j+2}}$ and $w_1 w_2 \cdots w_k$ is of length $j$, which in particular gives
$$w_3w_4\cdots w_j w_{j+1} \overline{w_{j+2}} =w_1\cdots w_j.$$
The result then follows by recalling $w_{j+1} \overline{w_{j+2}}=01$.
\end{proof}

%

Now we can prove Theorem \ref{isogen}.

\noindent{\em Proof of Theorem \ref{isogen}:} In this proof we will call the vertex labeled $*$ in $L^{\{11\}}_w$ by $0$ for convenience.
When $k=3$, there are only five possible choices of $w$: $000, 001,010,100$ and $101$. One can easily check that the corresponding unlabeled versions of $L_{w}^{\{11\}}$ for these five words are never isomorphic.

So we assume $k\geq 4$ in the following.
We will prove the following equivalent statement: given any (directed) graph
$G = (V, E)$ that arises as an unlabeled version of $L^{\{11\}}_w$ for some
$w \in \euscript{B}_k(\{11\})$, we can uniquely determine the word $w$, or equivalently, the
corresponding labeling of the graph $G$. 

Let $G = (V,E)$ be such a graph. Suppose we identified the vertex named $k-1$. Then we can use the induction argument in the beginning of the proof of Theorem \ref{theorem: isomorphism of graphs} to identify the vertices names $j$ for $1\leq j \leq k-1$. Finally by Fact \ref{g2w2}, we would have uniquely determined the word $w$. Thus all we have to do is to identify the vertex $v^*$ with name $k-1$. 
vertex $v^* \in V$ which must have label $k-1$.

For $i = 0, 1$, define
\begin{equation} V_i = \{v \in V: \text{outdegree}_G(v) = 2-i\}. \end{equation}
If there exists a vertex $v \notin V_0 \cup V_1$, i.e. if there is a vertex with out-degree zero, which must be the only such vertex by the construction of $L_w^{\{11\}}$, then by Fact \ref{g2w2} $v$ must have label $k-1$, as desired. So assume that $V_0, V_1$ are a partition of $V$. By Fact \ref{g2w2}, each vertex $x \in V_0$ must have vertex label $i(x)$ with $w_{i(x)} = 0$, and each vertex $y \in V_1$ --  except possibly the vertex with label $k-1$ -- must have vertex label $i(y)$ with $w_{i(y)} = 1$. Thus, since $w$ is an allowed word in $X^{\{11\}}$, each vertex $y \in V_1$ with $w_{i(y)} = 1$ must have $w_{i(y)+1} = 0$, so the outgoing edge at $y$ has endpoint $x \in V_0$. If there exists $v \in V_1$ such that the outgoing edge $v \to z$ has $z \notin V_0$, then $v$ must be the unique vertex with label $k-1$, and again we are done.

So assume that all outgoing edges from $V_1$ have endpoints in $V_0$ (and $V_0, V_1$ form a partition of $V$). Consider the following six cases:
\begin{align*}
\mbox{(Case 1): } w_1w_2w_3= 001;  \quad & \mbox{(Case 2): } w_1w_2w_3= 100;  & \mbox{(Case 3): } w_1 w_2 w_3 = 101; \\
\mbox{(Case 4): } w_1 w_2 w_3 = 000; \quad & \mbox{(Case 5): } w_1 w_2 w_3 w_4=0100; \quad & \mbox{(Case 6): } w_1 w_2 w_3 w_4=0101.
\end{align*}

One easily checks that, because $w \in \euscript{B}_k(\{11\})$, these are the only possible length-$3$ prefixes of $w$ (note that Case 5 and 6 combined corresponds to $w_1w_2w_3=010$).

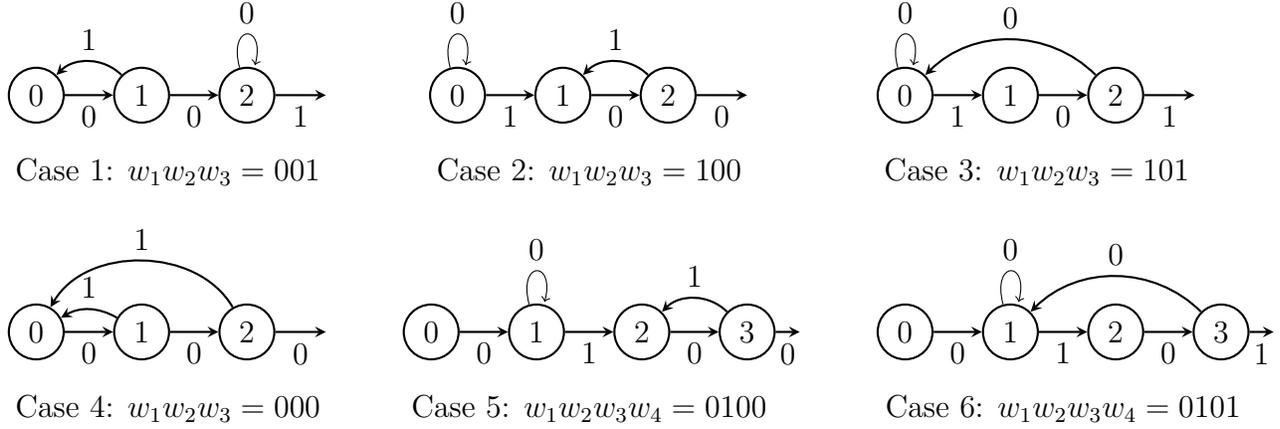
\begin{figure}[H]
\begin{center}
\begin{tikzpicture}[scale=0.35]
\draw [opacity=0.0] (0,5) grid (50,20);
\node [circle, draw, thick] (a1) at (0,19) {$0$};
\node [circle, draw, thick] (b1) at (4,19) {$1$};
\node [circle, draw, thick] (c1) at (8,19) {$2$};
\draw [-stealth, black, thick] (a1) -- node[below] {$0$} (b1);
\draw [-stealth, black, thick] (b1) -- node[below] {$0$} (c1);
\draw [-stealth, black, thick] (c1) -- node[below] {$1$} (11,19);
\draw [-stealth, black, thick] (b1) to [out=135, in=45] node[above] {$1$}(a1);
\path (c1) edge [loop above] node {$0$} (c1);
\node at (5, 17)[below] {Case 1: $w_1w_2w_3=001$} ;

\node [circle, draw, thick] (a2) at (16,19) {$0$};
\node [circle, draw, thick] (b2) at (20,19) {$1$};
\node [circle, draw, thick] (c2) at (24,19) {$2$};
\draw [-stealth, black, thick] (a2) -- node[below] {$1$} (b2);
\draw [-stealth, black, thick] (b2) -- node[below] {$0$} (c2);
\draw [-stealth, black, thick] (c2) -- node[below] {$0$} (27,19);
\path (a2) edge [loop above] node {$0$} (a2);
\draw [-stealth, black, thick] (c2) to [out=135, in=45] node[above] {$1$}(b2);
\node at (21, 17)[below] {Case 2: $w_1w_2w_3=100$} ;

\node [circle, draw, thick] (a3) at (33,19) {$0$};
\node [circle, draw, thick] (b3) at (37,19) {$1$};
\node [circle, draw, thick] (c3) at (41,19) {$2$};
\draw [-stealth, black, thick] (a3) -- node[below] {$1$} (b3);
\draw [-stealth, black, thick] (b3) -- node[below] {$0$} (c3);
\draw [-stealth, black, thick] (c3) -- node[below] {$1$} (44,19);
\path (a3) edge [loop above] node {$0$} (a3);
\draw [-stealth, black, thick] (c3) to [out=135, in=45] node[above] {$0$}(a3);
\node at (38, 17)[below] {Case 3: $w_1w_2w_3=101$} ;

\node [circle, draw, thick] (a4) at (0,10) {$0$};
\node [circle, draw, thick] (b4) at (4,10) {$1$};
\node [circle, draw, thick] (c4) at (8,10) {$2$};
\draw [-stealth, black, thick] (a4) -- node[below] {$0$} (b4);
\draw [-stealth, black, thick] (b4) -- node[below] {$0$} (c4);
\draw [-stealth, black, thick] (c4) -- node[below] {$0$} (11,10);
\draw [-stealth, black, thick] (b4) to [out=150, in=30] node[above] {$1$}(a4);
\draw [-stealth, black, thick] (c4) to [out=120, in=60] node[above] {$1$}(a4);
\node at (5, 8)[below] {Case 4: $w_1w_2w_3=000$} ;

\node [circle, draw, thick] (a5) at (15,10) {$0$};
\node [circle, draw, thick] (b5) at (19,10) {$1$};
\node [circle, draw, thick] (c5) at (23,10) {$2$};
\node [circle, draw, thick] (d5) at (27,10) {$3$};
\draw [-stealth, black, thick] (a5) -- node[below] {$0$} (b5);
\draw [-stealth, black, thick] (b5) -- node[below] {$1$} (c5);
\draw [-stealth, black, thick] (c5) -- node[below] {$0$} (d5);
\draw [-stealth, black, thick] (d5) -- node[below] {$0$} (29,10);
\path (b5) edge [loop above] node {$0$} (b5);
\draw [-stealth, black, thick] (d5) to [out=135, in=45] node[above] {$1$}(c5);
\node at (21, 8)[below] {Case 5: $w_1w_2w_3w_4=0100$} ;

\node [circle, draw, thick] (a6) at (33,10) {$0$};
\node [circle, draw, thick] (b6) at (37,10) {$1$};
\node [circle, draw, thick] (c6) at (41,10) {$2$};
\node [circle, draw, thick] (d6) at (45,10) {$3$};
\draw [-stealth, black, thick] (a6) -- node[below] {$0$} (b6);
\draw [-stealth, black, thick] (b6) -- node[below] {$1$} (c6);
\draw [-stealth, black, thick] (c6) -- node[below] {$0$} (d6);
\draw [-stealth, black, thick] (d6) -- node[below] {$1$} (47,10);
\path (b6) edge [loop above] node {$0$} (b6);
\draw [-stealth, black, thick] (d6) to [out=135, in=45] node[above] {$0$}(b6);
\node at (39, 8)[below] {Case 6: $w_1w_2w_3w_4=0101$} ;

\end{tikzpicture}
\end{center}
\caption{All six cases}
\label{all}
\end{figure}

Now, the subgraph (of $L_{w}^{\{11\}}$) we see at vertices $0,1,2$ (and vertices $0,1,2,3$ for cases (5.1) and (5.2)) are all given in Figure \ref{all}.

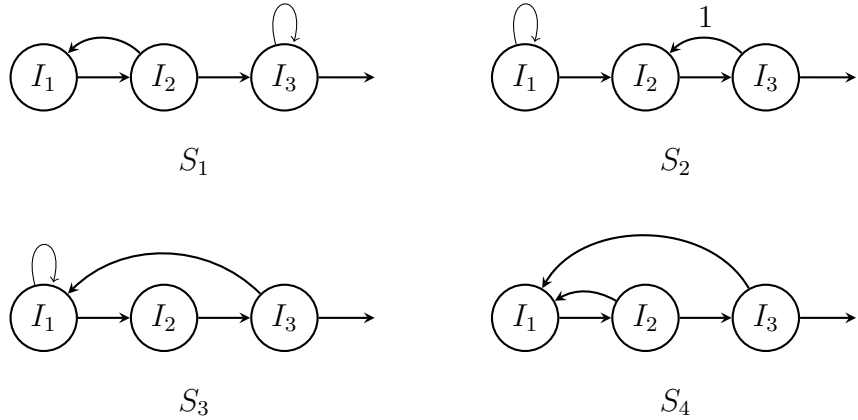
\begin{figure}[H]
\begin{center}
\begin{tikzpicture}[scale=0.4]
\draw [opacity=0] (0,7) grid (35,20);
\node [circle, draw, thick] (a1) at (3,19) {$I_1$};
\node [circle, draw, thick] (b1) at (7,19) {$I_2$};
\node [circle, draw, thick] (c1) at (11,19) {$I_3$};
\draw [-stealth, black, thick] (a1) -- node[below] {} (b1);
\draw [-stealth, black, thick] (b1) -- node[below] {} (c1);
\draw [-stealth, black, thick] (c1) -- node[below] {} (14,19);
\draw [-stealth, black, thick] (b1) to [out=135, in=45] node[above] {}(a1);
\path (c1) edge [loop above] node {} (c1);
\node at (8, 17)[below] {$S_1$} ;

\node [circle, draw, thick] (a2) at (19,19) {$I_1$};
\node [circle, draw, thick] (b2) at (23,19) {$I_2$};
\node [circle, draw, thick] (c2) at (27,19) {$I_3$};
\draw [-stealth, black, thick] (a2) -- node[below] {} (b2);
\draw [-stealth, black, thick] (b2) -- node[below] {} (c2);
\draw [-stealth, black, thick] (c2) -- node[below] {} (30,19);
\path (a2) edge [loop above] node {} (a2);
\draw [-stealth, black, thick] (c2) to [out=135, in=45] node[above] {$1$}(b2);
\node at (24, 17)[below] {$S_2$} ;

\node [circle, draw, thick] (a3) at (3,11) {$I_1$};
\node [circle, draw, thick] (b3) at (7,11) {$I_2$};
\node [circle, draw, thick] (c3) at (11,11) {$I_3$};
\draw [-stealth, black, thick] (a3) -- node[below] {} (b3);
\draw [-stealth, black, thick] (b3) -- node[below] {} (c3);
\draw [-stealth, black, thick] (c3) -- node[below] {} (14,11);
\path (a3) edge [loop above] node {} (a3);
\draw [-stealth, black, thick] (c3) to [out=135, in=45] node[above] {}(a3);
\node at (8, 9)[below] {$S_3$} ;

\node [circle, draw, thick] (a4) at (19,11) {$I_1$};
\node [circle, draw, thick] (b4) at (23,11) {$I_2$};
\node [circle, draw, thick] (c4) at (27,11) {$I_3$};
\draw [-stealth, black, thick] (a4) -- node[below] {} (b4);
\draw [-stealth, black, thick] (b4) -- node[below] {} (c4);
\draw [-stealth, black, thick] (c4) -- node[below] {} (30,11);
\draw [-stealth, black, thick] (b4) to [out=150, in=30] node[above] {}(a4);
\draw [-stealth, black, thick] (c4) to [out=120, in=60] node[above] {}(a4);
\node at (24, 9)[below] {$S_4$} ;

\end{tikzpicture}
\end{center}
\caption{``Base graphs" $S_1, S_2, S_3, S_4$}
\label{unlabeled}
\end{figure}

Let $S_1, S_2, S_3, S_4$ be the unlabeled subgraphs shown in Figure \ref{unlabeled},
which will be called the ``base graphs". Then, one can see from Figure \ref{all} that
for any $w$, at least one of the base graphs will appear in $L_{w}^{\{11\}}$. 
%

Now, given $G$, we build the following
algorithm to identify the vertex $0$, which in turn helps us identify the vertex $k-1$
(because $k-1$ is the terminal vertex of the only Hamiltonian path of
$L_{w}^{\{11\}}$ starting at the vertex $0$). The proof of the validity of the algorithm,
which is little tedious, can be found in the Appendix.

The algorithm is:

{\bf Step 1:}  First, we search for $S_1$ in $G$. If we see $S_1$ as a subgraph of $G$ then the vertex $I_1$ (in $S_1$) must be labeled $0$ (and then we are done).

{\bf Step 2:} If we can not find $S_1$ in $G$, then we search for $S_2$. If we see $S_2$ somewhere, we consider two subcases:
\begin{itemize}
\item if there is a vertex (call this vertex $T$) going into $I_1$ (of $S_2$) and $T$ has in-degree $0$, then we are in Case 5, and the vertex $T$ must be labeled $0$ and we are done;
\item otherwise, if all vertices going into $I_1$ have in-degree at least $1$, then we must be in Case 2, where the vertex with a self-loop must be labeled $0$.
\end{itemize}

{\bf Step 3:} If neither $S_1$ nor $S_2$ is found in $G$, then  we search for $S_3$. If we see $S_3$ somewhere, then again we consider two subcases:
\begin{itemize}
\item if there is a  vertex (call this vertex $T$) going into $I_1$ and $T$ has in-degree $0$, then we are in Case 6 and the vertex $T$ must be labeled $0$;
\item otherwise, if all vertices going into $I_1$ have in-degree at least $1$, then we must be in Case 3, and the vertex with a self-loop must be labeled $0$.
\end{itemize}

{\bf Step 4:} If we fail to find $S_1, S_2$ and also $S_3$ in $G$, then we must see $S_4$. In this case, $I_1$ (in $S_4$) must be labeled $0$ and again we are done. \qed

\section {Conjugacy of one-dimensional SFTs obtained by forbidding one word}\label{section:conjugacy}


In this section, we 
consider the relation between  $X_{\{u\}}$ and $X_{\{v\}}$ from a symbolic dynamics point of view. In particular, we seek conditions that imply that $X_{\{u\}}$ and $X_{\{v\}}$ are conjugate (conjugacy was defined in Section~\ref{section:intro}). Note that if there is a permutation on the alphabet $[q]$ which takes $u$ to $v$, then $X_{\{u\}}$ is trivially conjugate to $X_{\{v\}}$. In Proposition \ref{two_operations} we identify a large class of non-trivial examples.

One simple way to get a conjugacy between $X_{\{u\}}$ and $X_{\{v\}}$ is to find a $1$-block conjugacy with $1$-block inverse from $X_{L_{\{u\}}}$ to $X_{L_{\{v\}}}$. Unfortunately,  Corollary \ref{recover_word} tells us that 
conjugacies of this type do not exist if $u$ and $v$ are essentially different. However, a different type of conjugacy between $X_{\{u\}}$ and $X_{\{v\}}$ may exist. Indeed, our main result in this section shows that, if $u$ and $v$ both have only trivial self-overlap, then, other than the four exceptional cases of Proposition \ref{lem:irr}, $X_{\{u\}}$ is always conjugate to $X_{\{v\}}$. Moreover, it turns out that a version of this result also holds when the ambient shift space is the golden mean shift.

We need a definition of overlap between two different words. For two words $u, v\in [q]^k$, the {\em overlap set of $u$ and $v$} is defined by
$$
\ov(u,v)=\{i\in [1,k]: v_{[1,i]}=u_{[k-i+1, k]}\}.
$$
Note that when $u=v$, this becomes the self-overlap of $u$ defined in Section \ref{section:earlier_1d}.


\subsection{The ambient shift is the full shift} \label{conj_subsec1}

We start with the following proposition which says that if $u$ and $v$ both have trivial self-overlap and they do not overlap with each other, then $X_{\{u\}}$ and $X_{\{v\}}$ are conjugate to each other. 

\begin{pr} \label{two_operations}
	Let $u, v\in [q]^k$. 
If $$\mathrm{overlap}(u,u) = \mathrm{overlap}(v,v)={\{k\}}\quad \mbox{and} \quad \mathrm{overlap}(u,v) = \mathrm{overlap}(v,u) =\emptyset,$$
then $X_{\{u\}}$ is conjugate to $X_{\{v\}}$.
\end{pr}
\begin{proof}

Denote $u=u_1u_2\cdots u_k$ and $v=v_1v_2\cdots v_k$. Define a $(2k-1)$-block code $\Phi:[q]^{2k-1} \rightarrow [q]$ with memory $k-1$ and anticipation $k-1$
by
\begin{align}
\Phi(x_{[-k+1, k-1]})=
\begin{cases}
u_{1-i} \qquad &\mbox{if $x_{[i, i+k-1]}=v$ for some $i\in [-k+1,0]$} \\
v_{1-i} \qquad &\mbox{if $x_{[i, i+k-1]}=u$ for some $i\in [-k+1,0]$} \\
x_0 \qquad &\mbox{otherwise}.
\end{cases} \label{replace_coor}
\end{align}
Note that $\Phi$ is well-defined because $u, v$ both have trivial self-overlap and no cross-correlation.
Let $\phi$ be the sliding block code induced by $\Phi$. Thus, 
we deduce from (\ref{replace_coor}) that for all $i$,
\begin{align}
\phi(x)_{[i,i+k-1]}=
\begin{cases}
u \qquad \mbox{if $x_{[i, i+k-1]}=v$} \\
v \qquad \mbox{if $x_{[i, i+k-1]}=u$},
\end{cases} \label{replace_uv}
\end{align}
i.e., $\phi$ replaces appearances of $u$ in $x$ with $v$ and vice versa.

We first claim that
\begin{align} \label{replace_not_uv}
\phi(x)_{[i,i+k-1]} \notin \{u,v\} \mbox{ if } x_{[i,i+k-1]}\notin\{u,v\} \qquad \mbox{for all $i$.}
\end{align}
To see this, for any fixed $i$, consider two cases: if there exists $j\in [i-k+1, i+k-1]$ such that $x_{[j, j+k-1]}\in \{u,v\}$, then by (\ref{replace_uv}), $\phi(x)_{[j, j+k-1]}\in \{u,v\}$ and therefore $\phi(x)_{[i,i+k-1]}\notin \{u,v\}$ since $\mathrm{overlap}(u,u)=\mathrm{overlap}(v,v)=\{k\}$ and $\mathrm{overlap}(u,v)=\mathrm{overlap}(v,u)=\emptyset$; on the other hand, if for all $j\in [i-k+1, i+k-1]$, $x_{[j, j+k-1]}\notin \{u,v\}$, then we directly infer from (\ref{replace_coor}) that  $\phi(x)_{[i,i+k-1]}=x_{[i,i+k-1]}\notin \{u,v\}$.


Now, by (\ref{replace_uv}) and(\ref{replace_not_uv}), $\phi$ is an involution and therefore a conjugacy on $[q]^\mathbb{Z}$.  
It remains to show $\phi(X_{\{u\}})=X_{\{v\}}$. We first show $\phi(X_{\{u\}})\subset X_{\{v\}}$. To see this, let $x\in X_{\{u\}}$. We claim that $\phi(x)\in X_{\{v\}}$. If not, there exists $i$ such that $\phi(x)_{[i, i+k-1]}=v$. Then, according to (\ref{replace_uv}) and (\ref{replace_not_uv}), we must have $x_{[i,i+k-1]}=u$, which contradicts $x\in X_{\{u\}}$. A similar argument gives $\phi(X_{\{v\}})\subset X_{\{u\}}$, which is equivalent to $X_{\{v\}}\subset \phi(X_{\{u\}})$ since $\phi$ is an involution. Thus, $\phi(X_{\{u\}})=X_{\{v\}}$ and the proof is complete.
\end{proof}

\begin{rem}
The conjugacy we constructed in the proof is given by replacing any appearance of $v$ by $u$ and vice versa. It is not hard to check that this type of conjugacy between $X_{\{u\}}$ and $X_{\{v\}}$ extends to a conjugacy on the full $q$-shift. 
We also note that  if $u$ can be obtained from $v$ by permuting the alphabet, then there is another type of conjugacy between $X_{\{u\}}$ and $X_{\{v\}}$.
We call these two types of conjugacies together {\em swap conjugacies}.
\end{rem}


To present our main results,
we start with the binary alphabet case. Suppose $q=2$. Let $k$ be a positive integer and define
	$$
	\mathcal{C}(k):=\{{w}\in \{0,1\}^k \setminus \{10^{k-1}, 01^{k-1}, 1^{k-1}0, 0^{k-1}1\}: \mathrm{overlap}(w,w)=\{k\} 
	\}.
	$$
Note that here we have to exclude the four words $10^{k-1}, 01^{k-1}, 1^{k-1}0, 0^{k-1}1$ since according to Proposition \ref{lem:irr}, each of these words gives an reducible SFT which can not be conjugate to any irreducible SFT.

Our first main result in this section is
\begin{thm}\label{conj_trivial_self_over}
	Let $q=2$. Given any two distinct words $u,v\in \mathcal{C}(k)$, there exist an integer $N\leq 6$ and a sequence of words $u=u^{(1)}, u^{(2)},\cdots, u^{(N-1)}, u^{(N)}=v$ in $\mathcal{C}(k)$ such that for any $1\leq i< N$, either $u^{(i+1)}$ is obtained from $u^{(i)}$ by flipping all bits or
	$$
	\ov(u^{(i)}, u^{(i+1)})=\ov(u^{(i+1)}, u^{(i)})=\emptyset.
	$$
	Consequently, there is a chain of swap conjugacies from $X_{\{u\}}$ to $X_{\{v\}}$, and the length of the chain is upper bounded by $6$.
\end{thm}

We need a series of lemmas to prove Theorem \ref{conj_trivial_self_over}. To introduce these lemmas, we first define
$$
\mathcal{C}_1(k)=\{w\in \mathcal{C}(k): w_1=1, w_k=0\}.
$$
For each $2\leq j\leq k-2$ and $2\leq l \leq k-2$, we also define sets 
\begin{align*}
D(j)&:=\{w \in \mathcal{C}_1(k): w_{[1,\lceil k/2 \rceil +j]}=1^j 0^{ \lceil k/2 \rceil}\}, \\
E(l)&:=\{w\in \mathcal{C}_1(k): w_{[\lceil k/2 \rceil-l+1, k]} = 1^{ \lfloor k/2 \rfloor }0^l\}.
\end{align*}

Let $z:= 1010^{k-3}$ and for each $2\leq j\leq k-2$, let $p(j):=1^j0^{k-j}$ be the ``test words". Note that $z\in \mathcal{C}_1(k)$ if $k>4$ and $p(j)\in D(j)\subset \mathcal{C}_1(k)$ if $k\geq 4$ and $2\leq j\leq\lfloor k/2\rfloor$. The following fact can be verified immediately.
\begin{fact} \label{test_words}
Let $k>4$. Then, for any $2\leq j\leq k-2$,
\begin{align*}
&\ov(z,z)= \ov(p(j), p(j)) = \{k\} \\
\mbox{and} \qquad &\ov(p(j),z)=\ov(z,p(j))=\emptyset.
\end{align*}
\end{fact}

We then have the following lemma.

\begin{lem} \label{not_special_head}
Let $k>4$. Then, for any $w\in \mathcal{C}_1(k)$, 
we have
\begin{enumerate}
\item If $\ov(p(\lfloor k/2 \rfloor), w)\neq \emptyset$, then $w\in D(j)$ for some $2\leq j\leq \lfloor k/2 \rfloor$;
\item If $\ov(w, p(\lfloor k/2 \rfloor)) \neq \emptyset$, then $w\in E(l)$ for some $2\leq l\leq \lceil k/2 \rceil$.
\end{enumerate}
\end{lem}

\begin{proof}
We prove item 1 here and the proof of item 2 is similar.

Suppose $\ov(p(\lfloor k/2 \rfloor), w)\neq \emptyset$. Since $w$ starts with $1$ and ends with $0$, $p(\lfloor k/2 \rfloor)$ must overlap with $w$ in at least $\lceil k/2 \rceil+1$ many bits. On the one hand, this overlapping substring is a suffix of $p(\lfloor k/2 \rfloor)$ and therefore it is of the form $1^j0^{\lceil k/2 \rceil}$ for some $1\leq j\leq \lfloor k/2 \rfloor$;  on the other hand, this overlapping substring is also a prefix of $w$, indicating that $w$ must start with $1^j0^{\lceil k/2 \rceil}$ for some $1\leq j\leq \lfloor k/2 \rfloor$. Now to complete the proof, it remains to show $j\neq 1$. To see this, suppose to the contrary that $j=1$. Then $w$ starts with $10^{\lceil k/2 \rceil}$. Noting that $w$ has trivial self-overlap, we must have $w=10^{k-1}$, which is impossible.
\end{proof}

An immediate consequence of Lemma \ref{not_special_head} is the following.

\begin{co} \label{contrapositive_head}
Let $k>4$ and $w\in \mathcal{C}_1(k)$. If $w\notin D(j)$ for any $2\leq j \leq \lfloor k/2\rfloor$ and $w\notin E(l)$ for any $2\leq l \leq \lceil k/2 \rceil$, then
$$
\ov(p(\lfloor k/2 \rfloor),w)=\ov(w,p(\lfloor k/2 \rfloor))=\emptyset.
$$
\end{co}


Lemma \ref{special_head} and Fact \ref{special_tail} below deal with strings in $D(j)$ and $E(l)$, respectively.
\begin{lem} \label{special_head}
Let $k>4$. Then, for any $2\leq j\leq \lfloor k/2 \rfloor$ and any $w\in D(j)$,
either $w=p(j)$ or
$$
\ov(w, p(j))=\ov(p(j),w)=\emptyset.
$$
\end{lem}

\begin{proof}
Suppose $w\neq p(j)$. We will show $\ov(w, p(j))=\ov(p(j),w)=\emptyset$.

Now, $\ov(p(j), w) = \emptyset$ is straightforward.
To see $\ov(w, p(j))=\emptyset$, suppose to the contrary that $\ov(w, p(j))\neq \emptyset$. Recalling that $w_k=0$, we know $w$ must have a suffix $1^j 0^{l}$ for some $l\geq 1$; on the other hand,  $w\in D(j)$ implies that $1^j 0^{\lceil k/2 \rceil}$ is a prefix of $w$. Since $w\neq p(j)$, we must have $l<\lceil k/2 \rceil$. But then $1^j0^l$ is both a prefix and a suffix of $w$, contradicting the fact that $w$ has only trivial self-overlap.
\end{proof}

The following fact is also easy to verify.
\begin{fact} \label{special_tail}
Let $k>4$. Then, for any $2\leq l{\color{blue}\leq} \lceil k/2 \rceil$ and any $w\in E(l)$,
$$
\mathrm{overlap}(w, z)=\mathrm{overlap}(z,w)=\emptyset.
$$
\end{fact}

We are now ready to prove Theorem \ref{conj_trivial_self_over}.

\noindent{\bf Proof of Theorem \ref{conj_trivial_self_over}:}
We only prove the first claim of Theorem \ref{conj_trivial_self_over}, since the claim on the existence of the swap-conjugacy chain then follows from Proposition \ref{two_operations}.

We first note that $\mathcal{C}(k)=\emptyset$ if $1\leq k\leq 3$ and $\mathcal{C}(k)=\{0011,1100\}$ if $k=4$. Therefore, Theorem \ref{conj_trivial_self_over} trivially holds when $k\leq 4$.

Now consider $k>4$.
Let $p(j), z, \mathcal{C}_1(k),D(j), E(l)$ be define as before, where $2\leq j, l\leq k-2$. Let
$$
\mathcal{C}_2(k)=\{w\in \mathcal{C}(k): w_1=0, w_k=1\}.
$$
Note that $\mathcal{C}(k)=\mathcal{C}_1(k)\cup \mathcal{C}_2(k)$ and there is a one-to-one correspondence between words in $\mathcal{C}_1(k)$ and words in $\mathcal{C}_2(k)$ given by flipping all bits. 
Thus, it suffices to prove Theorem \ref{conj_trivial_self_over} for $u,v\in \mathcal{C}_1(k)$ and for $N\leq 5$.

To this end, let $u\in \mathcal{C}_1(k)$ be arbitrary. We first claim that there exist
$N\leq 3$ and a sequence of words $u=u^{(1)}, u^{(2)},\cdots, u^{(N)}=z$ such that for any $1\leq i< N$, $u^{(i)}\in \mathcal{C}_1(k)$ and $\ov(u^{(i)}, u^{(i+1)})=\ov(u^{(i+1)}, u^{(i)})=\emptyset$.
Indeed, if $u=p(j)$ for some $2\leq j\leq k-2$, then, according to Fact \ref{test_words}, the claim holds with $u^{(1)}=u, u^{(2)}=z$;
if $u\notin D(j)$ for any $2\leq j\leq\lfloor k/2 \rfloor$ and $u\notin E(l)$ for any $2\leq l \leq \lceil k/2 \rceil$, then, by Corollary \ref{contrapositive_head} and Fact \ref{test_words}, the claim holds with $u^{(1)}=u, u^{(2)}=p(\lfloor k/2 \rfloor), u^{(3)}=z$;
if $u\in D(j)$ for some $2\leq j\leq\lfloor k/2 \rfloor$ with $u\neq p(j)$, then, by Lemma \ref{special_head} and Fact \ref{test_words}, the claim holds with $u^{(1)}=u, u^{(2)}=p(j), u^{(3)}=z$;
if $u\in E(l)$ for some $2\leq l\leq\lceil k/2 \rceil$, then by Fact \ref{special_tail}, the claim again holds with $u^{(1)}=u, u^{(2)}=z$. Thus, the claim is proved.

Now, by applying the argument above to any $v\in \mathcal{C}_1(k)$, we deduce that there exist $M\leq 3$ and a sequence of words $z=v^{(1)}, v^{(2)},\cdots, v^{(M-1)}, v^{(M)}=v$ such that for any $1\leq i< M$, $v^{(i)}\in \mathcal{C}_1(k)$ and $\ov(v^{(i)}, v^{(i+1)})=\ov(v^{(i+1)}, v^{(i)})=\emptyset$.
The proof is then completed by combining the two sequences $u^{(1)},\cdots, u^{(N)}$ and $v^{(1)}, \cdots, v^{(M)}$.
 \qed

We now generalize the conjugacy result in Theorem \ref{conj_trivial_self_over} to larger alphabets. 

\begin{thm}\label{theorem: conjugacy}
Suppose $q>2$. Let $u, v\in [q]^k$ such that $\ov(u,u)=\ov(v,v)=\{k\}$. Then, there exist an integer $N\leq 11$ and a sequence of words $u=u^{(1)}, u^{(2)},\cdots, u^{(N-1)}, u^{(N)}=v$ with trivial self-overlap such that for any $1\leq i< N$, either $u^{(i+1)}$ is obtained from $u^{(i)}$ by permuting the alphabet or
	$$
	\ov(u^{(i)}, u^{(i+1)})=\ov(u^{(i+1)}, u^{(i)})=\emptyset.
	$$
	Consequently, there is a chain of swap conjugacies from $X_{\{u\}}$ to $X_{\{v\}}$, and the length of the chain is upper bounded by $11$.
\end{thm}

\begin{rem}
Compared to Theorem \ref{conj_trivial_self_over}, words of the form $a^{k-1}b$ and $ab^{k-1}$ are not excluded in Theorem \ref{theorem: conjugacy} because by Remark \ref{Remark: Reducible subshift} the SFTs given by them are irreducible (rather than reducible as in the $q=2$ case) when $q>2$.
\end{rem}

\begin{proof}
Again, we only prove the first part of the theorem, since the existence of the swap-conjugacy chain then follows from Proposition \ref{two_operations}.	
	
For any $a, b\in [q]$ with $a\neq b$, define
$$
\mathcal{C}_{ab}(k):=\{w\in [q]^k: w_1=a, w_k=b, \ov(w,w)=\{k\}\}.
$$
Observe that for any such pair $a$ and $b$, there is a one-to-one correspondence between words in $\mathcal{C}_{10}(k)$ and $\mathcal{C}_{ab}(k)$ by permuting the alphabet. Thus, it suffices to prove Theorem \ref{theorem: conjugacy} for $u, v\in \mathcal{C}_{10}(k)$ and for $N\leq 9$. We remind the reader that words in $C_{10}(k)$ are not necessarily binary words.

We first consider the case $k\leq 4$. Suppose $k=3$. To prove Theorem \ref{theorem: conjugacy} holds for any $u,v\in \mathcal{C}_{10}(k)$, it suffices to show that for any word $w\in \mathcal{C}_{10}(k)$, there is a sequence of words $w=w^{(1)}, w^{(2)},\cdots, w^{(N)}=1^{k-1}0$ with trivial self-overlap such that for any $1\leq i< N$, either $w^{(i+1)}$ is obtained from $w^{(i)}$ by permuting the alphabet or $\ov(w^{(i)}, w^{(i+1)})=\ov(w^{(i+1), w^{(i)}})=\emptyset$. Indeed, if $w$ uses all but the letter $c\in [q]$, then, by choosing $w^{(1)}=w, w^{(2)}=c^{k-1}0, u^{(3)}=1^{k-1}0$, we have
$
\ov(w^{(1)},w^{(2)})= \ov(w^{(2)}, w^{(1)})=\emptyset
$
and $w^{(3)}$ is obtained from $w^{(2)}$ by permuting the alphabet; if $w\in \mathcal{C}_{10}(k)$ uses all the letters from $[q]$, then $q=3$ and we must have $w=120$. Thus, choosing $w^{(1)}=w=120, w^{(2)}=110$, we have $\ov(w^{(1)},w^{(2)})=\ov(w^{(2)}, w^{(1)})=\emptyset$. This proves the theorem for $k=3$.
The case $k=4$ can be handled similarly.


Now consider $k>4$. Let $u\in \mathcal{C}_{10}(k)$ be arbitrary. Let $z:= 1010^{k-3}$ which in particular is in $\mathcal{C}_{10}(k)$. We claim that, there exist $N\leq 5$ and a sequence of words $u=u^{(1)}, u^{(2)},\cdots, u^{(N-1)}, u^{(N)}=z$ with trivial self-overlap  such that for any $1\leq i< N$, either $u^{(i+1)}$ is obtained from $u^{(i)}$ by permuting the alphabet or $\ov(u^{(i)}, u^{(i+1)})=\ov(u^{(i+1)}, u^{(i)})=\emptyset$.

To prove this claim,  we consider three cases:
\begin{itemize}
\item If $u$ is binary and $u\notin \{10^{k-1}, 01^{k-1}, 1^{k-1}0, 0^{k-1}1\}$, then by the proof of Theorem \ref{conj_trivial_self_over}, the claim holds with $N=3$.

\item If $u\in \{10^{k-1}, 01^{k-1}, 1^{k-1}0, 0^{k-1}1\}$, we assume without loss of generality that $u=10^{k-1}$. Let $a\in [q]\setminus \{0,1\}$ and choose $u^{(1)}=u=10^{k-1}, u^{(2)}=aa0^{k-2}, u^{(3)}=110^{k-2}, u^{(4)}=z$ (note that these words all have trivial self-overlap). Then, one can check that $\ov(u^{(1)},u^{(2)})=\ov(u^{(2)}, u^{(1)})=\emptyset$, $u^{(3)}$ is obtained from $u^{(2)}$ by permuting the alphabet, and $\ov(u^{(3)},u^{(4)})=\ov(u^{(4)}, u^{(3)})=\emptyset$. Thus, the claim holds with $N=4$.

\item if $u$ contains a symbol in $[q]\setminus\{0,1\}$ 
, let $x(u)=1^{l+1}0^{k-l-1}$ where $l$ is the longest all-$1$ substring appearing in $u$. Then, $x(u)\in \mathcal{C}_{10}(k)$, $x(u)$ is binary and one check that $\ov(x(u), u)=\ov(u,x(u))=\emptyset$. Moreover, according to the two cases we just proved above, there exist a sequence of words $x(u)=u^{(1)}, u^{(2)},\cdots, u^{(K)}=z$ with trivial self-overlap such that $K\leq 4$ and for any $1\leq i< K$, either $u^{(i+1)}$ is obtained from $u^{(i)}$ by permuting the alphabet or $\ov(u^{(i)}, u^{(i+1)})=\ov(u^{(i+1)}, u^{(i)})=\emptyset$. Therefore, by adding $u$ to this sequence of words, we conclude that the same thing holds for the new sequence $u, u^{(1)},\cdots, u^{(K)}$, i.e., the claim holds with $N=5$.
\end{itemize}

Now that the claim is proved, we can apply it to any $v\in \mathcal{C}_{10}(k)$ to deduce that there exist $M\leq 5$ and a sequence of words $z=v^{(1)}, v^{(2)},\cdots, v^{(M-1)}, v^{(M)}=v$ with trivial self-overlap such that for any $1\leq i< M$, either $v^{(i+1)}$ is obtained from $v^{(i)}$ by permuting the alphabet or $\ov(v^{(i)}, v^{(i+1)})=\ov(v^{(i+1)}, v^{(i)})=\emptyset$.
The proof is then completed by combining the two sequences $u^{(1)},\cdots, u^{(N)}$ and $v^{(1)}, \cdots, v^{(M)}$.
\end{proof}


Finally, it is natural to ask whether we can extend Proposition \ref{two_operations} and Theorem \ref{conj_trivial_self_over} to the case when $u$ and $v$ have nontrivial self-overlap. Indeed, the extension of the former is possible while the extension of the latter is not.

\begin{pr} \label{two_operations_exten}(An extension of Proposition \ref{two_operations})
	Let $u, v\in [q]^k$  such that
	$$
	\mathrm{overlap}(u,u)=\mathrm{overlap}(v,v).
	$$
		If $$
		\mathrm{overlap}(u,v)=\mathrm{overlap}(v,u)= \mathrm{overlap}(u,u) \setminus \{k\},
		$$
		then there is a swap conjugacy between $X_{\{u\}}$ and $X_{\{v\}}$.
\end{pr}

\begin{proof}
The proof is similar to that of Proposition {\ref{two_operations}}. Indeed, one can check from the proof there that
$\mathrm{overlap}(u,v)=\mathrm{overlap}(v,u)= \mathrm{overlap}(u,u) \setminus \{k\}$ is sufficient to define a swap conjugacy which replaces any appearance of $u$ by $v$ and vice versa.
\end{proof}

The extension of Theorem \ref{conj_trivial_self_over} to words with non-trivial self-overlap is however not true, as can be seen from the following example.
\begin{exmp} \label{reverse_exam}
	Suppose $q=2$. Let $k=6$ and $\mathrm{overlap}(u,u)=\mathrm{overlap}(v,v)=\{3,6\}$. 
	Observe that $110110, 011011, 100100, 001001$ are the only four binary strings of length 6 having self-overlap set $\{3,6\}$. If we take, for example, $u=110110$ and $v=011011$, then one can check that no matter how you choose the intermediate words in between $u$ and $v$ (so that they form a chain of words from $u$ to $v$), there always exist two adjacent words in this chain such that Proposition \ref{two_operations_exten} fails. Therefore there does not exist a chain of swap conjugacies between $X_{\{u\}}$ and $X_{\{v\}}$.

\end{exmp}

\subsection{The ambient shift is the golden mean shift} \label{conj_subsec2}

	Results in the previous subsection can be extended to the case when the ambient shift is the golden mean shift by a similar argument. We present them in this subsection.


Let $X_{\{11\}}$ denote the golden mean shift, i.e., the set of binary bi-infinite sequences that do not contain $11$. Let $\euscript{B}_k({\{11\}})$ be the set of allowed words of length $k$ in $X_{\{11\}}$. For any $w\in \euscript{B}_k(\{11\})$ let $(X_{\{11\}})_{\{w\}}$ be the SFT obtained by forbidding $w$ from the golden mean shift.

The following proposition is an analog of Proposition \ref{two_operations} for the golden mean ambient shift case.

\begin{pr} \label{two_operations_2}
Let $u, v\in \euscript{B}_{k}({\{11\}})$ be such that

\begin{enumerate}
\item They have the same initial symbol and terminal symbol, i.e., $u_1=v_1$ and $u_k=v_k$;
\item They both have trivial self-overlap and they do not overlap with each other, that is,
$$\mathrm{overlap}(u,u) = \mathrm{overlap}(v,v)=\{k\}\quad \mbox{and} \quad \mathrm{overlap}(u,v) = \mathrm{overlap}(v,u) =\emptyset.$$
\end{enumerate}
Then there is a swap conjugacy between $(X_{\{11\}})_{\{u\}}$ and $(X_{\{11\}})_{\{v\}}$.
\end{pr}

\bigskip
We omit the proof since it is similar to the proof of Proposition \ref{two_operations}. Item 1 of Proposition \ref{two_operations_2} exactly means $u$ and $v$ have the same ``extender set" (see, for example, Section 2 of \cite{garcia2019extender} for the definition of extender sets). This condition is introduced here to ensure the point obtained from the replacement operation is in the golden mean shift as long as the original point is. We also point out that though the extender set condition is crucial for our techniques, it does not seem to be necessary for the conjugacy between SFTs obtained by forbidding one word (from a general ambient SFT).

We introduce some notations before generalizing Theorem \ref{conj_trivial_self_over}.
Let
\begin{align*}
\mathcal{R}:=\{w\in \euscript{B}_k({\{11\}}): w\in \{10^{k-1}, 0^{k-1}1 , (10)^{\frac{k-1}{2}}0, 0(01)^{\frac{k-1}{2}}\}\}
\end{align*}
where we regard $(10)^{\frac{k-1}{2}}0$ and $0(01)^{\frac{k-1}{2}}$ as empty words when $k$ is even. {One directly checks from $L_{w}^{\{11\}}$ that $(X_{\{11\}})_{\{w\}}$ is reducible when $w\in \mathcal{R}$ and $k\geq 3$, and the SFTs $\{(X_{\{11\}})_{\{w\}}: w\in \mathcal{R}\}$ are never conjugate to each other.} 
Define
\begin{align*}
\mathcal{G}(k)&:=\{{w}\in \euscript{B}_k(\{11\}): \mathrm{overlap}(w,w)=\{k\} \mbox{ and } w\notin \mathcal{R}\};\\
\mathcal{G}_1(k)&:=\{w\in \mathcal{G}(k): w_1=1, w_k=0\}, \qquad \mathcal{G}_2(k):=\{w\in \mathcal{G}(k): w_1=0, w_k=1\}.
\end{align*}

Our main theorem in this subsection is the following, which is an analog of Theorem \ref{conj_trivial_self_over}.

\begin{thm} \label{conj_GM}
For any $u, v\in \mathcal{G}_1(k)$ (resp. $u,v\in \mathcal{G}_2(k)$), there exist an integer $N\leq 5$ and a sequence of words $u=u^{(1)}, u^{(2)},\cdots, u^{(N-1)}, u^{(N)}=v$ in $\mathcal{G}_1(k)$ (resp. $\mathcal{G}_2(k)$) such that for any $1\leq i< N$, 
	$$
	\ov(u^{(i)}, u^{(i+1)})=\ov(u^{(i+1)}, u^{(i)})=\emptyset.
	$$
	Consequently, there is a chain of swap conjugacies between $(X_{\{11\}})_{\{u\}}$ and $(X_{\{11\}})_{\{v\}}$, and the length of the chain is upper bounded by $5$.
\end{thm}

We give a sketch of the proof of Theorem \ref{conj_GM} since it is similar to that of Theorem \ref{conj_trivial_self_over}. The main difference between these two proofs is the choice of ``test words".

\begin{proof}

We only prove the case when $u,v\in \mathcal{G}_1(k)$ since there is a one-to-one correspondence between elements in  $\mathcal{G}_1(k)$ and elements in $\mathcal{G}_2(k)$ given by reversing the words, and for two words $u$ and $v$, condition 2 in Proposition \ref{two_operations_2} is preserved after reversing them.

We assume $k\geq 7$ in the following since one can readily check that $\mathcal{G}_1(k)$ contains at most one element when $k<7$. We also assume $k$ is even to avoid notational cumbersomeness. \footnote{The case $k$ is odd, just as in the proof of Theorem \ref{conj_trivial_self_over}, can be handled by a similarly argument.}

Let $z:= 10010^{k-4}$ and define the ``test words" $p(j):=(10)^j 0^{k-2j}$ for all $2\leq j\leq \frac{k-2}{2}$. Note that 
these words are all in $\mathcal{G}_1(k)$. Since there are at least three consecutive $0$'s at the end of each $p(j)$, we have for any $2\leq j\leq \frac{k-2}{2}$,
$$
\ov(z,p(j))=\ov(p(j),z)=\emptyset.
$$

Now let $m:= \lceil\frac{k}{4} \rceil$ and define sets
\begin{align*}
D(j)&:= \{w\in \mathcal{G}_1(k): w_{[1,k+2j-2m]}=(10)^j 0^{k-2m}\} \qquad \mbox{for all $2\leq j\leq m$} \\
E(l)&:=\{w\in \mathcal{G}_1(k): w_{[k-2m-l+1,k]}=(10)^m 0^l\} \qquad \mbox{for all $2\leq l \leq k-2m$}.
\end{align*}
Then,  by a similar argument as in the proof of Lemma \ref{not_special_head}, we have for any $w\in \mathcal{G}_1(k)$, 
the following hold:
\begin{enumerate}
\item If $\mathrm{overlap}(p(m), w)\neq \emptyset$, then $w\in D(j)$ for some  $2\leq j\leq m$;
\item If $\mathrm{overlap}(w, p(m))\neq \emptyset$, then $w\in E(l)$ for some  $2\leq l \leq k-2m$.
\end{enumerate}
Therefore, for any $w\in \mathcal{G}_1(k)$ such that $w\notin D(j)$ for any $2\leq j \leq m$ and $w\notin E(l)$ for any $2\leq l \leq k-2m$, we must have
$$
\ov(w, p(m))=\ov(p(m),w)=\emptyset.
$$

We then deal with words in $D(j)$ and $E(l)$. Indeed, these two cases can be handled similarly as in Lemma \ref{special_head} and Fact \ref{special_tail}: if $w\in D(j)$ for some $2\leq j\leq m$, then a similar proof to Lemma \ref{special_head} gives that either $w=p(j)$ or $\mathrm{overlap}(w, p(j))=\mathrm{overlap}(p(j),w)=\emptyset$; 
if $w\in E(l)$ for some $2\leq l\leq k-2m$, then, similar to Fact \ref{special_tail}, we have  $\mathrm{overlap}(w, z)=\mathrm{overlap}(z,w)=\emptyset$.

Combining all the cases above, we conclude that for any $w\in \mathcal{G}_1(k)$, there is a sequence of words $w=w^{(1)}, \cdots, w^{(N)}=z$ in $\mathcal{G}_1(k)$ such that $N\leq 3$ and for any $1\leq i< N$, $\ov(w^{(i)}, w^{(i+1)})=\ov(w^{(i+1)}, w^{(i)})=\emptyset$, from which the theorem follows. 
\end{proof}

\begin{rem}
It is worth noting that in Theorem \ref{conj_GM} we need to consider the two disjoint sets $\mathcal{G}_1(k)$ and $\mathcal{G}_2(k)$ separately. This is because when the ambient shift is the golden mean shift, words in $\mathcal{G}_1(k)$ can not be obtained from words in $\mathcal{G}_2(k)$ by simply flipping the bits (while this is always possible when the ambient is the full $2$-shift). As a simple example, let $k=5$. Then $u=10100$ is in $\mathcal{G}_1(k)$ but the word $v=01011$, obtained by flipping all the bits of $u$, is not in $\mathcal{G}_2(k)$ because it is not allowed in the golden mean shift.
\end{rem}

For $k\geq7$, choose $w=10010^{k-4} \in \mathcal{G}_1(k)$. One can check from the graph $L_w^{\{11\}}$ that $(X_{\{11\}})_{\{w\}}$ is an irreducible SFT. Then, according to Theorem \ref{conj_GM}, $(X_{\{11\}})_{\{u\}}$ is irreducible for all $u\in \mathcal{G}_1(k)$. A similar argument also applies to $\mathcal{G}_2(k)$. Thus, we have the following corollary.

\begin{co}
Let $w\in \euscript{B}_k (\{11\})$ be such that $\ov(w,w)=\{k\}$. Then, $(X_{\{11\}})_{\{w\}}$ is reducible if and only if $w\in \mathcal{R}$.
\end{co}

\section{Higher dimensional shift spaces}
\label{section:higherd}
In this section, we will have a look at $\Z^d$ for $d\geq 2$. For patterns on a higher dimensional lattice it is already a question what the set of overlaps mean and what kind of overlap sets can arise. Further there is no nice analogue of Perron-Frobenius theory which will work in our context. Consequently everything becomes significantly more complicated and it is very rare to have general results for higher dimensional actions. This makes the results of this section, though special,  significant.

Given two sets $S, T\subset \Z^d$ and patterns $u\in [q]^S, u'\in [q]^T$, we will say that $u$ and $u'$ \emph{agree} if for all $i\in S\cap T$, we have that $u_i=u'_i$. Observe that patterns with disjoint support agree because they don't overlap anywhere.

Given a pattern $u$ on a set $S\subset \Z^d$, its translate the pattern $\sigma^{-i}(u)$ has a chance to non-trivially agree with $u$ only if their supports, the sets $S+i$ and $S$ intersect. If $j$ is a point in the intersection then $j=k+i$ for some $k\in S$. Rewriting this we find that $i=j-k\in S-S$.  Thus trivially we have that the translate of a shape $S$ by a vector $i$ intersects itself if and only if $i\in S-S$. In addition if the pattern and its translate by $i$ agree then we must have $u_j=u_{k}=u_{j-i}=(\sigma^{-i}(u))_j$. With this in mind we have the following definition of self-overlap set in higher dimension. To distinguish from the $d=1$ case we will call them agreement/self-agreement sets thereafter.
\begin{definition}
	Given a finite set $S\subset \Z^d$ and a pattern $u\in [q]^S$, we define
	$$\iv(u)=\{i\in S-S~:~ \text{ for all $j\in S\cap (S+i)$ for which $u_j=u_{j-i}$}\}.$$
	In other words, $\iv(u)$ is the set of sites $i\in S-S$ for which $u$ and its translate $\sigma^{-i}(u)$ agree on $S\cap(S+i)$.
\end{definition}

\begin{rem}
	\begin{enumerate}
		\item Recall the definition of the overlap sets introduced in Section \ref{section:earlier_1d} for words $w\in A^{\{1,2,\ldots, k\}}$. It is easy to see that
		$$\iv(w)=\{0\}\cup(k-\ov(w,w)) \cup(\ov(w,w)-k).$$
		In $d=1$ it was better to consider overlap sets because it is more intuitive to think of the matching suffixes and prefixes and because of its relationship with the correlation polynomial. In higher dimensions, we think about agreement sets instead because there is no natural way of talking about prefixes and suffixes.
		\item There is some redundancy in the definition: Indeed $i\in \iv(u)$ if and only if $-i\in \iv(u)$. But beyond this, we do not have a very deep understanding of what $\iv(u)$ could look like for a given set $S$ (see Section \ref{Section: weird overlap sets}). When $S$ is an interval in $\Z$ then they have been characterised by Guibas and Oldyzko \cite{guibas1981periods}.
	\end{enumerate}
\end{rem}

Suppose that we are given two patterns $u^{(0)}$ and $u^{(1)}$ where we want to replace all the appearances of $u^{(0)}$ by $u^{(1)}$. This does not necessarily make sense. For instance suppose $u^{(0)}=101$ and $u^{(1)}= 000$ and that we want to replace the appearances of $u^{(1)}$ by $u^{(0)}$ in the word  $w=1000010101$. Now $u^{(1)}$ appears starting from the second and third letter but we cannot concurrently place $u^{(0)}$ starting at the second and third position since the third letter can't concurrently be both $0$ and $1$. On the other hand, we can concurrently replace appearances of $u^{(0)}$ by $u^{(1)}$. This motivates the following proposition.

\begin{prop}\label{prop: replace when overlap contain}
Let $S, T\subset \Z^d$ and $u^{(0)}, u^{(1)}\in [q]^S$ be such that $\iv(u^{(0)})\subset \iv(u^{(1)})$. Suppose we are given a pattern $w\in [q]^T$ and $I\subset \Z^d$ be such that $w|_{i+S}=\sigma^{-i}(u^{(0)})$ for all $i\in I$. Then there exists $w'\in [q]^T$ such that $w'|_{i+S}=\sigma^{-i}(u^{(1)})$ for all $i\in S$ and $w'|_{T\setminus \cup_{i\in I}(i+S)}=w|_{T\setminus \cup_{i\in I}(i+S)}$.
\end{prop}

\begin{proof} For this proposition we need to check that the definition of $w'$ prescribed in the proposition makes sense. In particular for all $i_1, i_2\in I$ we need to check that the patterns $\sigma^{-i_1}(u^{(1)})$ and $\sigma^{-i_2}(u^{(1)})$ agree on $(i_1+S)\cap (i_2+S)$. This is trivially true if $(i_1+S)\cap (i_2+S)=\emptyset$. Suppose instead that $j\in (i_1+S)\cap( i_2+S)$. We need to verify that $(u^{(1)})_{j-i_1}=(u^{(1)})_{j-i_2}$.  Now since $w|_{i_1+S}=\sigma^{-i_1}(u^{(0)})$ and $\sigma^{-i_2}(u^{(0)})=w|_{i_2+S}$ and $(i_1+S)\cap(i_2+S)\neq \emptyset$ it follows that $i_2-i_1\in \iv(u^{(0)})$. By our hypothesis, $i_2-i_1\in \iv(u^{(1)})$ as well. Thus  $(u^{(1)})_{j-i_1}=(u^{(1)})_{j-i_1-i_2+i_1}=(u^{(1)})_{j-i_2}$. This completes the proof.
\end{proof}

Here on we will not be very careful about the support of the patterns. In particular when we say that $u\in A^S$ and $w|_{i+S}=u$, we just mean that
$w|_{i+S}=\sigma^{-i}(u)$. This is just to simplify the notation.

As before we will be interested in the properties of the shift space
$$X_{\{u\}}:=\{x\in [q]^{\Z^d}~:~u \text{ does not appear in }x\}.$$

Let us strengthen this assumption of containment a bit further so that it may apply to a variety of contexts where the ambient shift space is not the full shift. Let $X_\mathcal{F} \subset [q]^{\Z^d}$ be a shift space with a forbidden list $\mathcal F$. Given a finite set $S\subset \Z^d$ consider the set of \emph{globally allowed patterns} in $X_{\mathcal F}$ given by
$$\B_S(X_{\mathcal F})=\{x|_{S}~:~x\in X_{\mathcal F}\}$$
Let $(X_{\mathcal F})_{\{u\}}$ be obtained by forbidding a pattern $u\in \B_S(X_\mathcal{F})$:

$$(X_{\mathcal F})_{\{u\}}:=\{x\in X_{\mathcal F}~:~u \text{ does not appear in }x\}.$$

Let $u^{(0)}, u^{(1)}\in \B_S(X_{\mathcal F})$. We say that $u^{(0)}$ can be \emph{replaced} by $u^{(1)}$ in $X_{\mathcal F}$ if for all $x\in X_{\mathcal F}$ and a set of sites $I\subset \Z^d$ for which
$x|_{i+S}=u^{(0)}$ for all $i \in I$, there is a configuration $y\in X_{\mathcal F}$ such that
$$y|_{i+S}=u^{(1)} \text{ for all }i\in I\text{ and }y|_{\Z^d\setminus \cup_{i\in I}(i+S)}=x|_{\Z^d\setminus \cup_{i\in I}(i+S)}.$$

Let us make some elementary observations.

\begin{enumerate}
	\item We claim that a pattern $u^{(0)}$ can be replaced by a pattern $u^{(1)}$ in a full shift if and only if $\iv(u^{(0)})\subset \iv(u^{(1)})$. To see this, first observe that if  $\iv(u^{(0)})\subset \iv(u^{(1)})$ then $u^{(0)}$ can be replaced by $u^{(1)}$ in the full shift by by Proposition
\ref{prop: replace when overlap contain}. For the other direction suppose that $i\in \iv(u^{(0)})\setminus \iv(u^{(1)})$. Take any configuration $x\in [q]^{\Z^d}$ such that $x|_{S}=x|_{i+S}=u^{(0)}$. There exists no $y\in [q]^{\Z^d}$ such that $y|_{S}=y_{i+S}=u^{(1)}$ since $i \notin \iv(u^{(1)})$. Thus it follows that $u^{(0)}$ cannot be replaced by $u^{(1)}$ in the full shift.
	\item The interpretation of this replacement condition already becomes subtle if $X_{\mathcal F}$ is a non-trivial shift of finite type. Let us take a simple example. Suppose the alphabet is $[3]$ and consider $\mathcal F=\{0101,1010\}$, $u^{(1)}=101$ and $u^{(2)}=122$. Though $\iv(u^{(1)})\not\subset \iv(u^{(2)})$ (in fact $\iv(u^{(2)})\subset \iv(u^{(1)})$) we have that $u^{(1)}$ is replaceable by $u^{(2)}$. Let us see why. Let $x\in X_{\mathcal F}$ and $I \subset \Z$ such that $x|_{i+S}=u^{(1)}$ for all $i \in I$. First notice that if $i_1, i_2\in I$ then $(i_1+S)\cap (i_2+S)=\emptyset$. This is because any non-trivial agreement forms the pattern $10101$ which is forbidden in the shift space. This implies that replacing $u^{(1)}$ on the set $I$ by $u^{(2)}$ forms a valid point $y\in[3]^\Z$. Further since it only replaces some $0$'s and $1$'s by $2$'s it cannot possibly form a forbidden word which wasn't already there (which are composed of $0$'s and $1$'s). Thus $y\in X_{\mathcal{F}}$.
\end{enumerate}

How and when can patterns $u^{(0)}$ and $u^{(1)}$ be replaced to give a bijection is an interesting question. A little harder is to understand, when can we ensure that these replacements lead to conjugacy? This is related not only to self-agreements of the patterns but also agreements of the pattern $u^{(0)}$, $u^{(1)}$ and the ambient shift. We have discussed some of these issues in Section \ref{section:conjugacy} (look at Propositions \ref{two_operations} and \ref{two_operations_2}).

We will be interested in patterns $u^{(0)}$ and $u^{(1)}$ which are \emph{mutually replaceable }in $X$, meaning that $u^{(0)}$ can be replaced by $u^{(1)}$ and $u^{(1)}$ can be replaced by $u^{(0)}$ in $X$.

Many similar properties can be found in the literature. For instance Meyerovitch had a notion of exchangeability in \cite{meyerovitch2013gibbs} which allows for the replacement of a pattern $u$ by a pattern $v$ and vice versa one at a time rather than simultaneously and used it to prove that for all measures of maximal entropy $\mu$, $\mu([u])=\mu([v])$. Meyerovitch's condition is weaker and thus his conclusions follows for us as well. However our results do not hold under his weaker assumptions. Indeed, according to his definition, (say) the patterns $11$ and $10$ would be exchangeable in the full shift but $h(X_{\{11\}})>h(X_{\{10\}})$ by Proposition \ref{inequality_1d}. A similar property was also defined by Garcia-Ramos and Pavlov which allowed for replacements of a pattern $u$ by a pattern $v$ (and not vice versa) one at a time \cite{garcia2019extender} in points $x\in X$. However in our case, their agreement sets are different and that is why they are not mutually replaceable and by Proposition \ref{inequality_1d} we have that  $|\euscript{B}_{[1,n]}(X_{\{11\}})|>|\euscript{B}_{[1,n]}(X_{\{10\}})|$ for large enough $n$ (in fact all $n\geq 3$).

We will compare the size of the languages of various shift spaces in this section. Given a finite set $S\subset \Z^d$ and a shift space $X_{\mathcal F}$ and a pattern $u\in \B_U(X_{\mathcal F})$ for some $U$, we write
$$\B_S(X_{\mathcal F}, u)=\{x|_{S}~:~x\in X_{\mathcal F}\text{ and }u\text{ does not appear in }x|_S\}.$$
Note that if $\mathcal{F}$ is the empty set, $S$ is the interval $[1,n]$ and $u$ is a word, then $\B_S(X_{\mathcal{F}}, u)$ is just $\B_n(u)$ defined for the 1-dimensional case.

Now we are prepared to state the main results of this section.

\begin{thm}\label{Theorem: replace language equal}
Let $X_{\mathcal F}$ be a shift space and $u^{(0)},u^{(1)}$ be mutually replaceable in $X_\mathcal F$.  Then for all finite sets $S\subset \Z^d$ there is a bijection between
$\B_S (X_{\mathcal F}, u^{(0)})$ and $\B_S (X_{\mathcal F}, u^{(1)})$.
\end{thm}
\begin{rem}
\begin{enumerate}
	\item This generalises the implication ($1\Longrightarrow 6$) in Proposition~\ref{equality}.
	\item We believe that the same argument works for countable groups as well however we preferred to keep the simpler context which already illustrates the key issues.
\end{enumerate}
\end{rem}
\begin{proof}[Proof of Theorem~\ref{Theorem: replace language equal}]

The proof is a simple application of the inclusion-exclusion principle. Assume that $u^{(0)}, u^{(1)}\in [q]^{U}$ and that $0\in U$. Given a subset $T\subset S$ such that $T+U\subset S$ and $k=0,1$ consider
$$\B_S^{T,k}(X_{\mathcal F})=\{w\in \B_S(X_\mathcal F)~:~ w|_{i+U}= u^{(k)}\text{ for }i \in T\}.$$
Now consider the map $\phi_{T, k}: \B_S^{T,k}(X_{\mathcal F})\to \B_S^{T,1-k}(X_{\mathcal F})$ obtained by replacing appearances of $u^{(k)}$ by $u^{(1-k)}$, namely,
$$(\phi_{T,k}(w))|_{i+U} = u^{(1-k)}\text{ for $i\in T$ and }(\phi_{T,k}(w))|_{S\setminus \cup_{i\in T}(i+U)}=w|_{S\setminus \cup_{i\in T}(i+U)}.$$
This gives an injective map from $\B_S^{T,k}(X_{\mathcal F})$ to $\B_S^{T,1-k}(X_{\mathcal F})$ for $k=0,1$ proving that the cardinalities of the two sets are the same.

Now notice that the set $\B_S (X_{\mathcal F}, u^{(k)})$ is the set of patterns from $X_{\mathcal F}$ on $S$ in which $u^{(k)}$ does not appear even once. By inclusion-exclusion we have that
$$|\B_S (X_{\mathcal F}, u^{(k)})|=\sum_{T+U\subset S}(-1)^{|T|} |\B_S^{T,k}(X_{\mathcal F})|.$$
We remark that the summation on the right hand side also includes $T=\emptyset$, in which case, $\B_S^{T,k}(X_{\mathcal F}) = \B_S(X_{\mathcal F})$. Since $|\B_S^{T,0}(X_{\mathcal F})|=|\B_S^{T,1}(X_{\mathcal F})|$ for all $T\subset S$, this completes the proof.
%

\end{proof}

\begin{rem}
	One may be tempted to conjecture here that if $u^{(0)}$ is replaceable by $u^{(1)}$ in $X_{\mathcal F}$ then there is an injection from $\B_S (X_{\mathcal F}, u^{(0)})$ and $\B_S (X_{\mathcal F}, u^{(1)})$. This is not clear! Indeed we only get an injection from $\B_S^{T,0}(X_{\mathcal F})$ to $\B_S^{T,1}(X_{\mathcal F})$ and we cannot seem to use inclusion-exclusion directly.
	\end{rem}

\begin{rem}
	Note that the map does not give a direct bijection between $\B_S (X_{\mathcal F}, u^{(0)})$ and $\B_S (X_{\mathcal F}, u^{(1)})$ when $u^{(0)}$ is replaceable by $u^{(1)}$. We have observed a similar phenomenon in Section \ref{section:conjugacy}. More directly one can observe this when one tries to replace appearance of $u^{(0)}=1100$ by $u^{(1)}=1000$ in the full shift. Though both of them have the same (trivial) agreement set, notice that if one is to carry out the replacement in $11100$ then one obtains $11000$ which still has the word $u^{(0)}$. One can obviously apply the same map again but it is not clear if it would lead to a bijection. The reason is that the two patterns $u^{(0)}$ and $u^{(1)}$ have non-trivial agreement with each other (see Proposition \ref{two_operations}).
\end{rem}
The following corollary will follow easily from some fundamental properties of entropy.

\begin{co}\label{corollary: entropy matching}
Let $X_{\mathcal F}$ be a shift space and $u^{(0)},u^{(1)}$ be mutually replaceable in $X_\mathcal F$.  Then
$$h((X_{\mathcal{F}})_{\{u^{(0)}\}})=h((X_{\mathcal {F}})_{\{u^{(1)}\}}).$$
\end{co}
The proof will make use of the following concept. Given a forbidden list of patterns $\mathcal G$, the set of {\em locally allowed} patterns in $X_{\mathcal G}$ on the set $S$ is given by
$$\B_S^{loc}(\mathcal G):=\{a\in A^S~:~\text{ patterns from $\mathcal G$ do not appear in $a$}\}.$$
While this list of patterns is dependent on the specific choice of forbidden patterns $\mathcal G$, it follows that unlike the elements of globally allowed patterns, that is, $\B(X_{\mathcal G})$, the set of locally allowed patterns is algorithmically computable if $\mathcal G$ is finite.

\begin{proof}
	We have that
	$$\B_S(X_{\mathcal F, u^{(k)}})\subset \B_S(X_{\mathcal F}, u^{(k)})\subset \B_S^{loc}(\mathcal F\cup\{u^{(k)}\}).$$
It follows from \cite[Equation 3.9]{MR2129258},\cite{friedland1997entropy} that
\begin{eqnarray*}
h((X_{\mathcal F})_{\{u^{(k)}\}})&=&\lim_{n\to \infty}\frac{1}{|[1,n]^d|}\log\left(|\B_{[1,n]^d}(X_{\mathcal F, u^{(k)}})|\right)\\
&=& \lim_{n\to \infty}\frac{1}{|[1,n]^d|}\log\left(|\B_{[1,n]^d}^{loc}(\mathcal F\cup\{u^{(k)}\})|\right)
\end{eqnarray*}
Thus it follows that
$$h((X_{\mathcal F})_{\{u^{(k)}\}})=\lim_{n\to \infty}\frac{1}{|[1,n]^d|}\log\left(|\B_{[1,n]^d} (X_{\mathcal F}, u^{(k)})|\right).$$
The corollary now follows from Theorem \ref{Theorem: replace language equal}.
\end{proof}
A similar statement also follows for periodic points. Given a subgroup $\Lambda\subset \Z^d$, we denote the set of periodic points with period $\Lambda$ in a shift space $X_{\mathcal F}$ by
$$per^{\Lambda}(X_{\mathcal F})=\{x\in X_{\mathcal F}~:~\sigma^{i}(x)=x\text{ for all }i\in \Lambda\}$$
and for all $S\subset \Z^d$ its language
$$per^{\Lambda}_S(X_{\mathcal F})=\{x|_S~:~x\in X_{\mathcal F}, \sigma^{i}(x)=x\text{ for all }i\in \Lambda\}.$$
For instance if a sequence $x\in A^{\Z}$ has period $10$, the the corresponding subgroup is $\Lambda=10\Z$. Let $D\subset \Z^d$ be a fundamental domain for $\Lambda$, that is, for all $i \in \Z^d$ there exists a unique $\lambda\in \Lambda$ and $d\in D$ such that $d+\lambda=i$. For instance in the case $\Lambda=2\Z\times \{0\}$, a fundamental domain corresponding to it is $D=\{0,1\}\times \Z$. Let $X_{\mathcal F}$ be a shift space and $u \in \B_U(X_{\mathcal F})$ be a pattern for some finite $U \subset \Z^d$.

Define for all $S\subset D$, the restricted language of periodic points in which $u$ does not appear on S. More symbolically, we denote it by
$$per_S^{\Lambda}(X_{\mathcal F}, u)=\{x|_S~:~x\in per^{\Lambda}(X_{\mathcal F})\text{ and }x_{i+U}\neq u \text{ for all $i\in S$ }\}.$$

Note the subtle difference between definition of $per_S^{\Lambda}(X_{\mathcal F}, u)$ and $\B_S(X_{\mathcal F}, u)$. In the former we look at patterns on $S+U$ while in the latter we look at patterns on $S$. This is simply because in the former we are restricting to $S\subset D$ and even if $u$ appears in a periodic point $x\in per^{\Lambda}(X_{\mathcal F})$ it may not appear in $x|_{D}$. The following proposition generalises Theorem \ref{Theorem: replace language equal}.
\begin{prop}\label{proposition: replacement periodic points}
	Let $\Lambda\subset \Z^d$ be a subgroup and $D\subset \Z^d$ be a fundamental domain for $\Lambda$. Let $u^{(0)}$ and $u^{(1)}$ be patterns in a shift space $X_{\mathcal F}$ which are mutually replaceable. Then for all finite sets $S\subset D$ there is a bijection between $per^\Lambda_S(X_{\mathcal F}, u^{(0)})$ and $per_S^{\Lambda}(X_{\mathcal F}, u^{(1)})$.
\end{prop}
The proof is very similar to the proof of Theorem \ref{Theorem: replace language equal} and thus we only sketch the idea. For the set $T\subset S$ consider the set
$$per_S^{\Lambda, T, k}(X_{\mathcal F})=\{x|_S~:~x\in per^{\Lambda}(X_{\mathcal F})\text{ and }x|_{i+U}=u^{(k)}\text{ for }i\in T\}.$$
By using a map which replaces $u^{(0)}$ by $u^{(1)}$ on $i+\lambda+U$ for all $i \in T, \lambda\in \Lambda$ and similarly another one which replaces $u^{(1)}$ by $u^{(0)}$ we have that there is a bijection between the sets
$per_S^{\Lambda, T, 0}(X_{\mathcal F})$ and $per_S^{\Lambda, T, 1}(X_{\mathcal F})$. Finally by inclusion-exclusion we have a bijection between the sets $per_S^{\Lambda}(X_{\mathcal F}, u^{(0)})$ and $per_S^{\Lambda}(X_{\mathcal F}, u^{(1)})$.

Since $\Z^d$ acts on $per^{\Lambda}(X_{\mathcal F})$ with stabiliser $\Lambda$ we have a natural action of the group $\Z^d/\Lambda$ on $per^{\Lambda}(X_{\mathcal F})$. It follows that the topological entropy of the action of $\Z^d/\Lambda$ is simply equal to

$$h_{\Z^d/\Lambda}(X_{\mathcal F})=\lim_{n\to \infty}\frac{\log(|per^{\Lambda}_{ [1,n]^d}(X_{\mathcal F})|)}{\text{number of cosets of }\Lambda\text{ intersecting }[1,n]^d}.$$

 The following is an easy and direct corollary.

\begin{co}\label{corollary: periodic points}
	Let $\Lambda\subset \Z^d$ be a subgroup. Let $u^{(0)}$ and $u^{(1)}$ be patterns in a shift space $X_{\mathcal F}$ which are mutually replaceable. Then the topological entropy of $per^{\Lambda}((X_{\mathcal F})_{\{u^{(0)}\}})$ equals the topological entropy of $per^{\Lambda}((X_{\mathcal F})_{\{u^{(1)}\}})$ for the action of the group $\Z^d/\Lambda$. In addition the cardinality of $per^{\Lambda}((X_{\mathcal F})_{\{u^{(0)}\}})$ equals the cardinality of $per^{\Lambda}((X_{\mathcal F})_{\{u^{(1)}\}})$.
\end{co}
\begin{remark}
This generalises the implication $1\Longrightarrow 5$ from Proposition~\ref{equality} and generalises a consequence of a theorem of Lind~\cite[Theorem 1]{Lin89}.
\end{remark}
\begin{proof}[Proof of Corollary \ref{corollary: periodic points}]
	The proof of  the entropy statement is analogous to the proof of Corollary \ref{corollary: entropy matching}. By Proposition \ref{proposition: replacement periodic points} for all $S\subset D$ there is a bijection between $per^\Lambda_S(X_{\mathcal F}, u^{(0)})$ and $per_S^{\Lambda}(X_{\mathcal F}, u^{(1)})$.  Now as in Corollary  \ref{corollary: entropy matching} we can use ideas from \cite[Equation 3.9]{MR2129258},\cite{friedland1997entropy} to prove that the topological entropies must be equal. Finally note that
	as $S\uparrow D$,
	$$|per_S^{\Lambda}(X_{\mathcal F}, u^{(k)})|\uparrow |per^{\Lambda}((X_{\mathcal F})_{\{u^{(k)}\}})|$$
for both $k=1,2$. Given the bijection in Proposition \ref{proposition: replacement periodic points}, we see that the cardinalities must be equal. This completes the proof.
\end{proof}




\section{Discussion} \label{sec:discussion}

\subsection{The conjugacy problem for shift spaces with one forbidden word}
One of the major question which we weren't able to address in this paper is that of conjugacy.

\begin{ques}
Given two words $w, w' \in A^{k}$, is there an algorithm to determine whether the shift spaces $X_{\{w\}}$ and $X_{\{w'\}}$ are conjugate?
\end{ques}
It follows from Proposition \ref{equality} that a necessary condition is that the correlation polynomials for $w$ and $w'$ are identical. In the case that  $w$ and $w'$ have only trivial self-overlaps, Theorem \ref{theorem: conjugacy} provides a conjugacy. Beyond that, we don't know much. For instance, 
though we know from Example \ref{reverse_exam} that there does not exist a chain of swap conjugacies between the shift spaces $X_{\{110110\}}$ and $X_{\{011011\}}$, we don't know whether they are conjugate. In this case, the correlation polynomials are identical and one can check that the shift spaces agree on the main known conjugacy invariants (see~\cite[Section 12.3, Exercise 12.3.8]{LM21}).

	
\subsection{From follower set graph to correlation polynomial}

One significant question left open by our results is the implicit connection between the correlation polynomial $\phi_w$ and the graphs $L_w$, or more concretely, the $d_i(a)$ edge configurations. On one hand, the correlation polynomial is a natural object in terms of both hitting times and entropies. Namely, the explicit formula $\mathbb{E} \tau_w = q \phi_w(q)$ shows that each self-overlap of $w$ contributes to the expected hitting time of $w$ is an explicit way; and the connection between $\phi_w$ and the entropy $h(X_{\{w\}})$ is also concrete, if slightly less explicit, the entropy being the largest real root less than $q$ of the polynomial $1+(t-q)\phi_w(t)$. On the other hand, the condition $D(w,w')$ is also natural when comparing hitting times and entropies of two candidate shifts $X_{\{w\}}$ and $X_{\{w'\}}$, as Propositions \ref{pr:wine} and \ref{pr:cheese} show. These theorems are also more intuitive than the corresponding facts for the correlation polynomial: $D(w,w')$ guarantees that `more edges go further back in the chain $L$,' which heuristically suggests the Markov chain should take longer to hit the terminal state, and also that the graph should have a larger space of allowable paths (and thus higher entropy).

As a result, one expects to be able to convert between these two viewpoints. Indeed, by Theorem \ref{theorem: isomorphism of graphs}, each graph $L_w$ corresponds to a unique word $w$, and thus a unique correlation polynomial $\phi_w$. However, the proof of Theorem \ref{theorem: isomorphism of graphs} constructs this correspondence implicitly. Let $\mathcal{P}$ denote the family of polynomials that occur as correlation polynomials for some alphabet $[q]$ and word $w$ over $[q]$, and let $\mathcal{L}$ denote the set of graphs arising as follower set graphs of some such $w$. The family $\mathcal{P}$ can be described by explicit (but complex) `forward/backward propagation' rules, or alternatively by a recursive algorithm ~\cite{guibas1981periods}. Is it possible to give a similar characterization of the family $\mathcal{L}$? One attempt is to observe that the family $\mathcal{L}$ is closed under taking prefixes, i.e. if $L \in \mathcal{L}$, then the subgraph of $L$ obtained by deleting the vertex with label $k-1$, if $L = L_w$ for $w \in [q]^k$ is also in $\mathcal{L}$. Thus, one could hope that there exists an algorithm to decide which extensions of a given $L \in \mathcal{L}$ are also in $\mathcal{L}$. Realizing this sketch would require finding additional properties of $\mathcal{L}$ beyond those in Proposition \ref{prop:donkey_kong}.

Another issue with condition $D(w, w')$ is that it induces an incomplete partial ordering on the set of words $w \in [q]^k$, in the sense that there are pairs $w, w'$ such that neither $D(w, w')$ nor $D(w', w)$ holds, but $h(X_{\{w\}}) \neq h(X_{\{w'\}})$ ($w = 1001$ and $w' = 1100$ is an example). Additionally, the resulting poset is not `closed from below:' there exist words $w \in [q]^k$ whose entropy is strictly larger than the minimum possible entropy among words of length $k$, but there is no word $w'$ such that $D(w, w')$ holds. (A minimal example for the binary alphabet is $w = 1011$.) The poset is `closed from above,' in that $D(1^k, w)$ holds for all words $w \in [q]^k$. Is it possible to weaken the condition $D(w,w')$ to something which is compatible with follower set graphs, and which induces the same ordering on words $w \in [q]^k$ as ordering by entropy of the shifts $X_{\{w\}}$, or equivalently, the ordering of the hitting times $\tau_w$ by stochastic domination? Perhaps the alternative assumptions appearing in the main theorems of ~\cite{DS16} are a good starting point.



\subsection{When the ambient shift space is an SFT with memory one}
{ When the ambient shift space is the full $q$-shift, given words $u$ and $v$ of the same length, $h(X_{\{u\}})$ and $h(X_{\{v\}})$ can be easily compared using the correlation polynomial: $h(X_{\{u\}}) > h(X_{\{v\}})$ iff  $\phi_u(q) > \phi_v(q)$; see Proposition \ref{inequality_1d}. It is natural to ask if this result can be extended to some cases where the ambient shift space $X$ is not the full shift, by replacing $q$ with $\lambda$ where $h(X) = \log(\lambda)$. Numerical evidence suggests that at least a weak inequality of this type may be true (under suitable conditions). We show that this is actually true in some very special cases below (Proposition~\ref{eigen_geq_2_re1}).}

Let $X_{\mathcal{F}}$ be an irreducible  SFT with memory one over $[q]$. Let $T$ be the associated matrix of $X_{\mathcal{F}}$, i.e., it is a $q$-by-$q$ matrix and $T_{ij}=1$ if and only if $ij$ is allowed in $X_\mathcal{F}$. Let $\chi_T(t)$ denote the characteristic polynomial of $T$ and $\lambda_T$ be the largest root of $\chi_T(t)=0$.  For an allowed word $w$ in $X_T$ of length $k$, let $\phi_w(t)$ denote the correlation polynomial of the $2$-block representation of $w$ (the $2$-block representation of $w=w_1w_2\cdots w_k$ is $w^{[2]}=(w_1w_2) (w_2w_3)\cdots (w_{k-1} w_k)$, which is a $(k-1)$-block over the alphabet $[q]\times [q]$). The {\em extender set} of $w$ is the set of words that can be combined with $w$ to form a new allowed word in $X_T$.

 Let $(X_{\mathcal{F}})_{\{w\}}$ be the SFT obtained from $X_T$ by forbidding $w$.
 Just as in the case of the full shift, there is an analogous conjugacy presentation of  $(X_{\mathcal{F}})_{\{w\}}$, meaning a labeled graph s.t. the set of all bi-infinite label sequences is $(X_{\mathcal{F}})_{\{w\}}$ and the map from bi-infinite edge sequences to bi-infinite label sequences is one-to-one.  Moreover, there is an analogous
 characteristic polynomial $\chi_w(t)$, with nonzero constant term, whose largest root $\lambda_w$ satisfies $h(X_{\mathcal{F}})_{\{w\}}) = \log \lambda_w$ (for any conjugacy presentation).


Letting $i,j$ denote the initial, terminal symbols of $w$, one easily derives from
\cite[Theorem 1]{Lin89} that
\begin{equation}
\label{Doug}
\chi_w(t)= \chi_T(t)\phi_w(t) + cof_{ij}(tI-T)
\end{equation}
where $cof_{ij}(t-IT)$ equals $(-1)^{i+j}$ times the determinant of the matrix $t-IT$ with its $i$th row and $j$th column deleted, and

\begin{pr} \label{eigen_geq_2_re1}
	Let $X_\mathcal{F}$ be an ambient irreducible SFT with memory one, and $u, v$ be allowed words in $X_\mathcal{F}$ with length $k$ and the same extender set. Let $\eta_T$ be the second largest real root of $\chi_T(t)=0$. Then any of the following implies $\phi_u(\lambda_T)>\phi_v(\lambda_T)\iff h((X_{\mathcal{F}})_{\{u\}}) > h((X_{\mathcal{F}})_{\{v\}}) $:
	\begin{enumerate}
		\item[(a)] $\lambda_T\geq 2$, $\eta_T \geq 2$ and $\mbox{cof}_{ij} (\eta_T I -T)<0$; 
		\item[(b)] $\lambda_T>2$, $\eta_T\leq 2$ and $\mbox{cof}_{ij} (2I-T)<0$; (or $\lambda_T>2, \eta_T<2$ and $\mbox{cof}_{ij} (2I-T)\leq 0$)
		\item[(c)] $\lambda_T>2$, $\eta_T<2$, $\mbox{cof}_{ij}(2I-T)>0$ and $k\geq \log\frac{-\mbox{cof}_{ij}(2I-T)}{\chi_T(2)}+1$.
		\item[(d)] $\lambda_T\geq 2$ and $k$ is large enough.
	\end{enumerate}
\end{pr}

\begin{proof}
	It suffices to show $\phi_u(\lambda_T)>\phi_v(\lambda_T)\iff \lambda_u>\lambda_v $. Note that we have $\lambda_T\geq 2$ in all the four cases.
	We first write
	$$
	\phi_u(t)=\sum_{i=0}^{k-1} a_i t^i, \quad \phi_v(t)=\sum_{i=0}^{k-1} b_i t^i
	$$
	where $a_i, b_i \in \{0,1\}$. Define $r:= \max\{i: a_i=1, b_i=0\}$. This $r$ is well-defined since $\lambda_T>2$. Now $\phi_u(t)-\phi_v(t)>t^r -\sum_{i=0}^{r-1}t^i$. Let $\lambda^*$ be the largest root of $t^r -\sum_{i=0}^{r-1}t^i=0.$
	
$(\Rightarrow)$	: To prove this direction, we first claim it suffices to show $\lambda_v>\max\{\eta, \lambda^*\}$. For if $\lambda_v>\max\{\eta_T, \lambda_r\}$, then $\chi_T(\lambda_v)<0$ and $\phi_u(\lambda_v)-\phi_v(\lambda_v)>(\lambda_v)^r-\sum_{i=0}^{r-1}(\lambda_v)^i>0$. Thus, by (\ref{Doug}),
$$
\chi_u(\lambda_v)=\chi_u(\lambda_v) - \chi_v(\lambda_v) =\chi_T(\lambda_v)(\phi_u(\lambda_v) - \phi_v(\lambda_v)) < 0
$$
and therefore $\lambda_u>\lambda_v$.

Now, it is not hard to show that all the four condition implies $\lambda_v>\{\eta_T, \lambda_r\}$:

\noindent - Case {\em (a)}: Recall from (\ref{Doug}) that $\chi_v(\eta_T)=\chi_T(\eta_T) \phi_v(\eta_T)+\mbox{cof}_{ij}(\eta I-T)$. Now the conditions in $(a)$ imply $\chi_v(\eta_T)<0$, which gives $\lambda_v>\eta_T$. Observing that $\lambda^*<2$, we indeed have $\lambda_A>\max\{\eta_T, \lambda^*\}$.

\noindent - Case {\em (b)}: By a similar argument, it is not hard to derive $\chi_v(2)<0$. So $\lambda_v>2$ and therefore $\lambda_v>\max\{\eta_T, \lambda^*\}$.

\noindent - Case {\em (c)}: 
Note that $\chi_T(2)<0$ in this case. Since $\phi_v(2)\geq 2^{k-1}$, we derive from (\ref{Doug}) that $\chi_v(2)<\chi_T(2)2^{k-1}+\mbox{cof}_{ij}(2I-T)\leq 0$ because $k\geq \log\frac{-\mbox{cof}_{ij}(2I-T)}{\chi_T(2)}+1$ by assumption. Therefore, $\lambda_v>2$ which immediately gives $\lambda_v>\max\{\eta_T,\lambda^*\}$.

\noindent - Case {\em (d)}:
Irreducibility of $X$ implies $\eta_T<\lambda_T$. Also $\lambda^*<2\leq \lambda_T$. When $k$ is large enough, the entropy of $(X_{\mathcal{F}})_{\{v\}}$ is close to the entropy of $X_\mathcal{F}$, i.e., $\lambda_v$ is close to $\lambda_T$, which implies $\lambda_v>\max\{\eta_T, \lambda^*\}$.

$(\Leftarrow)$: Let $\lambda_u>\lambda_v$. We assume to the contrary that $\phi_u(\lambda_T)\leq \phi_v(\lambda_T)$. If $\phi_u(\lambda_T)< \phi_v(\lambda_T)$, then the proof of the other direction implies $\lambda_u<\lambda_v$, a contradiction; if $\phi_u(\lambda_T)=\phi_v(\lambda_T)$, then $\phi_u(t)=\phi_v(t)$ for all $t$ and therefore (\ref{Doug}) tells us $\chi_u(t)=\chi_v(t)$ for all $t$ (where we use the fact that the term $cof_{ij}(tI-T)$ is the same for $u$ and $v$ due to the same extender set assumption). Thus, $\lambda_u=\lambda_v$, a contradiction.
\end{proof}



\subsection{Agreement sets in higher dimensions}\label{Section: weird overlap sets}
We would like to discuss the question: Is there a ``nice'' characterisation of agreement sets for patterns in $\Z^d$? To illustrate what we mean by ``nice'', in the case when $d=1$ and the underlying shape is an interval, patterns are just words. The agreement sets of words are well-understood \cite[Theorem 5.1]{guibas1981periods}. This theorem gives us two facts.
\begin{enumerate}
	\item A nice algorithmic way of generating agreement sets.
	\item Two conditions: Forward and backward propagation which characterise what subsets of $\Z$ are agreement sets.
\end{enumerate}
As far as our study goes such a characterisation seems to be missing when $d\geq 2$. For us, this question is critical because it is the first step to understand how the agreement set of a pattern influences the drop in entropy when it is forbidden.  It also seems to be very important in other contexts, for instance, in pattern matching \cite{amir2022multidimensional}. The second motivation lead to extensive enquiry on agreement sets starting from the papers: \cite{periodicity_soda_92, regnier1993unifying}. Nevertheless a characterisation of agreement sets seems to be open. A nice survey on this and related problems can be found in \cite{smith2017properties}.

We will end with a couple of remarks:

\begin{enumerate}
	\item In $d=1$, \cite{guibas1981periods} showed that the collection of agreement sets of words in the binary alphabet equals the collection for arbitrarily sized alphabet. It is not enough to look at the binary alphabet to characterise agreement sets in higher dimensions. Indeed a quick enumeration will show that there exists no $w\in \{0,1\}^{[1,2]^2}$ such that $\iv(w)=\{(0,0)\}$ but clearly the agreement set of the pattern $\begin{smallmatrix}
		1&2\\3&4
	\end{smallmatrix}\in \{1,2,3,4\}^{[1,2]^2}$ has only $\{(0,0)\}$ as its self-agreement.
	\item If the shape $S\subset \Z^d$ is not convex then it may also happen for $w\in A^S$ that $u\in \iv(w)$ but $2u\notin \iv(w)$. One can already see this in $d=1$ when the set $S$ is not an interval but here is a connected example in $2$ dimensions. Consider the pattern
	$$w= \begin{smallmatrix}
		1&1\\
		1& \\
		1&0.
	\end{smallmatrix}$$
	It is clear $(0,1)\in \iv(w)$ but $(0,2)\notin \iv(w)$. This example shows that forward propagation (see \cite{guibas1981periods} for the definition) is not true for the agreement sets for general shapes in $\Z^2$.
\end{enumerate}
\subsection{Higher dimensional SFTs}\label{subsection: higher d inequalities}
In Section \ref{section:higherd} we discussed how the equality of the agreement sets implies equality of entropy. However a much more interesting and fundamental question is to understand how the agreement sets can cause a drop in entropy.

For instance if we are given $w, w'\in A^{[-n,n]^d}$, is there an easy way to decide whether or not $h(X_{\{w\}})>h(X_{\{w'\}})$? Towards this direction here is a conjecture that we were unable to prove.

\begin{conj}
	Let $w, w'\in A^{[-n,n]^d}$ be such that $\iv(w)\subsetneq \iv(w')$. Then $h(X_{\{w\}})<h(X_{\{w'\}})$.
\end{conj}

The conjecture even if proved brings us only a tiny bit closer to understanding how the agreement sets influence entropy. The first step of this program is to first understand what agreement sets look like; this has already been mentioned before in Section \ref{Section: weird overlap sets}. A baby version of the more general question (which does not satisfy the hypothesis of the conjecture), is to determine which of the following has larger entropy $X_{\{q^{(1)}\}}$ or $X_{\{q^{(2)}\}}$ where
$$q^{(1)}={\begin{smallmatrix}
		1&0&1&0&1&0\\
		0&1&0&1&0&1\\
		1&0&1&0&1&0\\
		0&1&0&1&0&1\\
		1&0&1&0&1&0\\
		0&1&0&1&0&1
\end{smallmatrix}} \text{ and }q^{(2)}=\begin{smallmatrix}
1&0&0&1&0&0\\
0&0&0&0&0&0\\
0&0&0&0&0&0\\
1&0&0&1&0&0\\
0&0&0&0&0&0\\
0&0&0&0&0&0
\end{smallmatrix}.
$$
Intuitively, we have that $q^{(1)}$ has `smaller' agreements than $q^{(2)}$ thus as in $d=1$ we conjecture that $h(X_{\{q^{(1)}\}})>h(X_{\{q^{(2)}\}})$. The more general version of such a conjecture will require a deeper understanding of agreement sets.

To show how much more complicated obtaining strict inequalities of entropy in higher dimensions tends to be, we will now give a proof sketch in a simple case which we know how to handle. In the following we assume that the reader is familiar with the basics of higher dimensional symbolic dynamics. Let $p^{(1)}$ be an $n\times n$ pattern with a single $1$ and rest as $0$'s
while $p^{(2)}$ is an $n\times n$ pattern full of $0$'s.
\begin{claim}
	$h(X_{\{p^{(1)}\}})<h(X_{\{p^{(2)}\}})$.
\end{claim}

Let us prove the claim. For this we will use the\emph{ lexicographic order} on $\Z^2$: $(i,j)>(k,l)$ if $j>l$ or $j=l$ and $i>k$. The reverse of this order is called the \emph{reverse lexicographic order}.

We begin by constructing an injective map from $\B_{[1,m]^2}(X_{\{p^{(1)}\}})$ to $\
B_{[1,m]^2}(X_{\{p^{(2)}\}})$. Let $p$ be a pattern in $\B_{[1,m]^2}(X_{\{p^{(1)}\}})$. In lexicographic order of appearance, replace the appearances of $p^{(2)}$ by $p^{(1)}$ one by one in $p$. Observe
\begin{enumerate}
	\item The replacement cannot create another $p^{(1)}$ except where the replacement takes place.
	\item The replacement cannot create an appearance of $p^{(2)}$.
\end{enumerate}
One can prove that this procedure is reversible because one can now replace appearances of $p^{(1)}$ by $p^{(2)}$ in the reverse lexicographic order to get back the pattern $p$. By continuing the replacement until all appearances of $p^{(2)}$ have been replaced by $p^{(1)}$ we get a pattern $p'\in \B_{[1,m]^2}(X_{\{p^{(2)}\}})$. Since the procedure is reversible, we have that the map $p\to p'$ is injective and
$$|\B_{[1,m]^2}(X_{\{p^{(1)}\}})|\leq |\B_{[1,m]^2}(X_{\{p^{(2)}\}})|.$$
Now we will observe that certain patterns from $\B_{[1,n+1]^2}(X_{\{p^{(2)}\}})$ cannot appear in such a $p'$. Start with the $(n+1)\times (n+1)$ pattern full of $0$'s. In the reverse lexicographic order, one by one, replace the appearances of $p^{(2)}$ by $p^{(1)}$ to get  a pattern $q$. Observe that $q\in \B(X_{\{p^{(2)}\}})$. However $q$ cannot appear in a pattern $p'$ constructed as above. For the sake of contradiction suppose that there did appear such a $p'$.  The proof splits into three cases:
\begin{enumerate}
	\item \emph{The $1$ appearing in $p^{(1)}$ appears at the upper right corner: } In this case $q$ is the pattern with four appearances of $1$ clumped in the upper right corner. To obtain $p$ from $p'$, we can now replace the appearance of $p^{(1)}$ by $p^{(2)}$ in $p'$ in the reverse lexicographic order. Only one of the appearances of $1$ in $q$ will be replaced by $0$, and so the resulting pattern $p$ will still have a $p^{(1)}$. Thus $p\notin \B_{[1,m]^2}(X_{\{p^{(1)}\}})$ which leads to a contradiction.
	\item\emph{$1$ appears in top or right edge in $p^{(1)}$: } This is similar to the case above and here too we obtain $p \notin \B_{[1,m]^2}(X_{\{p^{(1)}\}})$.
	\item \emph{$1$ does not appear in top or right edge on $p^{(1)}$: } Suppose the $1$ appears at the $(i,j)$ position in $p^{(1)}$ where $i<n$ and $j<n$. In this case $q$ is the pattern supported on the $(n+1)\times (n+1)$ box (say $B$) which is full of $0$'s except at $(i+1, j+1)$ where it is $1$. One can easily see that going from $p'$ to $p$, the $q$ is replaced by a pattern full of $0$'s on $B$. Now to reconstruct $p'$ we would have to replace the appearance of $p^{(2)}$ by $p^{(1)}$ in the lexicographic order. But then it follows that there is $(i',j')\in B$ which is lexicographically smaller than $(i+1, j+1)$ such that $p'_{(i',j')}=1$, a contradiction. 
\end{enumerate}
Thus we find that in fact
$$|\B_{[1,m]^2}(X_{\{p^{(1)}\}})|\leq |\B_{[1,m]^2}((X_{\{p^{(2)}\}})_{\{q\}})|.$$
It follows that $h(X_{\{p^{(1)}\}})\leq h((X_{\{p^{(2)}\}})_{\{q\}})$. To prove the claim it is enough to prove that $h((X_{\{p^{(2)}\}})_{\{q\}})<  h(X_{\{p^{(2)}\}})$. This follows because $1$ is a safe symbol\footnote{$1$ is a safe symbol is a $\Z^d$ shift space $X$ if for all $x\in X$ and $k \in \Z^d$ the configuration $y$ given by
$$y_l=\begin{cases}
x_l&\text{ if }l\neq k\\
1&\text{ if }l=k
\end{cases}$$ is an element of $X$} in $X_{\{p^{(2)}\}}$ which implies that $X_{\{p^{(2)}\}}$ is entropy minimal, that is, forbidding any pattern from $X_{\{p^{(2)}\}}$ causes a loss in entropy \cite{quas2003entropy}. This completes the sketch.

We see that the argument here is of a very different nature from those in $d=1$. At the same time it is somewhat special for the patterns $p^{(1)}$ and $p^{(2)}$ which have been especially chosen to demonstrate the method. The only other result known to us in this context is by Pavlov~ \cite{pavlov2011perturbations} who gives upper and lower bounds in the drop in entropy when a single pattern is forbidden, depending on the shape of the pattern, in the case of strongly irreducible subshifts.

\bibliographystyle{alpha}
\bibliographystyle{abbrv}
\bibliography{biblio}

\section*{Appendix: Validity of the algorithm in the proof of Theorem \ref{isogen}}\label{appB}

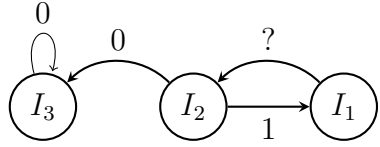
\begin{figure}[H]
\begin{center}
\begin{tikzpicture}[scale=0.5]
\draw [opacity=0] (0,18) grid (16,20);
\node [circle, draw, thick] (a1) at (3,19) {$I_3$};
\node [circle, draw, thick] (b1) at (7,19) {$I_2$};
\node [circle, draw, thick] (c1) at (11,19) {$I_1$};
\draw [-stealth, black, thick] (b1) -- node[below] {1} (c1);
\draw [-stealth, black, thick] (b1) to [out=135, in=45] node[above] {$0$}(a1);
\draw [-stealth, black, thick] (c1) to [out=135, in=45] node[above] {$?$}(b1);
\path (a1) edge [loop above] node {$0$} (a1);

\end{tikzpicture}
\end{center}
\caption{If $\mathcal{L}(I_1)=j+1, \mathcal{L}(I_2)=j, \mathcal{L}(I_3)=l$ for $S_1$, where $l<j$.}
\label{case1}
\end{figure}

{\em{Proof of Step 1}}: Suppose $S_1$ appears in $G$. 
We first show that the labels of $I_1, I_2, I_3$ must be $\mathcal{L}(I_1)=j, \mathcal{L}(I_2)=j+1, \mathcal{L}(I_3)=j+2$ for some $j\geq 0$ (here, for any vertex $J\in \mathcal{V}(G)$, $\mathcal{L}(J)$ denotes the label of $J$).  To see this, first observe that the labels of $I_1, I_2$ must be two adjacent integers since $I_1, I_2$ are in the same $2$-cycle. Now if $\mathcal{L}(I_1)=j+1, \mathcal{L}(I_2)=j$ for some $j\geq 0$, then the labeled version of $S_1$ is given in Figure \ref{case1} (where we use the fact that outgoing edges at $I_2$ have different labels). But (as you can see from the edge labeling in Figure \ref{case1}) the forward edge from $I_2$ to $I_1$ is labeled $1$, which, by Fact \ref{g2w2}, means $I_1$ should not have a ``going back edge'', a contradiction. Thus, we must have $\mathcal{L}(I_1)=j, \mathcal{L}(I_2)=j+1$. But then, the edge $I_2 \to I_3$ must be a ``forward edge" and therefore $\mathcal{L}(I_3)=j+2$.

Now we claim that $j=0$. 
If not, then there is a forward edge from $j-1$ to $j$. But then, Fact \ref{1cycle} implies this edge must be labeled $0$, and Fact \ref{2cycle} implies this edge must be labeled $1$, which is a contradiction.

Therefore, if we see $S_1$, then the vertex $I_1$ in $S_1$ (cf. Figure \ref{unlabeled}) must be labeled $0$.


{\em{Proof of Step 2}}: If  $S_2$ appears in $G$: 

Again, we first consider the label of the vertices $I_1, I_2, I_3$ in $S_2$. The existence of the $2$-cycle at $I_2, I_3$ implies that the labels of $I_2, I_3$ must be two adjacent integers. 
We claim that we must have $\mathcal{L}(I_2)<\mathcal{L}(I_3)$. If not, then the edge $I_3\to I_2$ must be a ``forward edge", which means $\mathcal{L}(I_3)=\mathcal{L}(I_2)-1$. On the other hand, noting that the edge $I_1\to I_2$ must also be a ``forward edge", we have $\mathcal{L}(I_1)=\mathcal{L}(I_2)-1$, which means $I_1$ and $I_3$ are the same vertex, a contradiction.


So assume $\mathcal{L}(I_2)=j+1, \mathcal{L}(I_3)=j+2$ for some $j$. Since $I_1$ has a self-loop, the edge from $I_1$ to $I_2$ must be a ``forward edge". Thus, we have $\mathcal{L}(I_1)=j$.

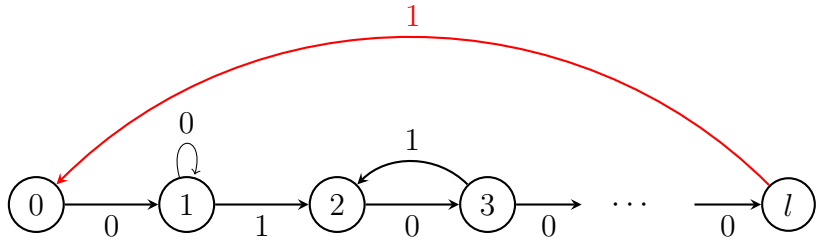
\begin{figure}[H]
\begin{center}
\begin{tikzpicture}[scale=0.5]
\draw [opacity=0] (0,17) grid (24,20);
\node [circle, draw, thick] (a1) at (1,19) {$0$};
\node [circle, draw, thick] (b1) at (5,19) {$1$};
\node [circle, draw, thick] (c1) at (9,19) {$2$};
\node [circle, draw, thick] (d1) at (13,19) {$3$};
\node at (16, 19)[right] {$\cdots$} ;
\node [circle, draw, thick] (e1) at (21,19) {$l$};
\draw [-stealth, black, thick] (a1) -- node[below] {$0$} (b1);
\draw [-stealth, black, thick] (b1) -- node[below] {$1$} (c1);
\draw [-stealth, black, thick] (c1) -- node[below] {$0$} (d1);
\draw [-stealth, black, thick] (d1) -- node[below] {$0$} (15.5,19);
\draw [-stealth, black, thick] (18.5,19) -- node[below] {$0$} (e1);
\draw [-stealth, black, thick] (d1) to [out=135, in=45] node[above] {$1$}(c1);
\draw [-stealth, red, thick] (e1) to [out=135, in=45] node[above] {$1$}(a1);
\path (b1) edge [loop above] node {$0$} (b1);

\end{tikzpicture}
\end{center}
\caption{An illustration of why the in-degree of the vertex $0$ is $0$ for Case 5. Here, the red edge is invalid.}
\label{case2graph}
\end{figure}

We then consider what value $j$ can take. By Fact \ref{1cycle} and Fact \ref{2cycle} , one have $j\in \{0,1\}$.
Indeed $j=0$ corresponds to Case 2 of Figure \ref{all}, $j=1$ corresponds to Case 5 of Figure \ref{all},
and the way to distinguish these two cases is to check the in-degree of vertices that can go to $I_1$: for case 2, any vertex that goes into the self-loop-vertex must have in-degree $\geq 1$; for case 5, the vertex $0$ goes into the self-loop-vertex and we claim that the vertex $0$ must have in-degree $0$. This is because if there is another vertex $l$ such that $d_l=0$, then $w_l\overline{w_{l+1}}=01$ so the ``going back edge'' at $l$ will never go to $0$ (it at least goes to the vertex $2$, see Figure \ref{case2graph}).


Having known how to distinguish Case 2 from Case 5, identifying the vertex labeled $0$ is straight forward: when we are in Case 2, then the vertex with a self-loop will be labeled $0$; if we are in Case 5, then there is a unique vertex with in-degree $0$ and there is an edge from this vertex to $I_1$; this vertex must be labeled $0$.

\medskip

{\em{Proof of Step 3}}: Suppose $S_1$ and $S_2$ do not appear in $G$, and $S_3$ appears. We first claim that the vertices of $S_3$ must be labeled $\mathcal{L}(I_1)=j, \mathcal{L}(I_2)=j+1, \mathcal{L}(I_3)=j+2$ for some $j\geq 0$. To see this, noting that $I_1$ has a self-loop, the edge $I_1 \to I_2$ must be a forward edge, which implies $\mathcal{L}(I_1)=j$ and $\mathcal{L}(I_2)=j+1$ for some $j\geq 0$. Now, if $\mathcal{L}(I_3)\neq j+2$, then the edge $I_2\to I_3$ must be a ``going back edge,'' from which, we infer that $\mathcal{L}(I_3)=j-1$. Now, by Fact \ref{1cycle}, the forward edge $I_1 \to I_2$ must be labeled $1$. But then, by Fact \ref{g2w2}, $I_2$ should not have a ``going back edge," a contradiction.


For the value of $j$, we first claim that $j\leq 2$. For if not, then the self-loop at the vertex $j$ implies $w_1 w_2 \cdots w_j=0^j$ where $j>2$. But then, one check that $w_{j+1}w_{j+2}\overline{w_{j+3}}=100$ and $w_{j-2} w_{j-1} w_j=000$, which means the ``going back edge" $j+2 \to j$ should not exist. 

\begin{figure}[H]
\begin{center}
\begin{tikzpicture}[scale=0.5]
\draw [opacity=0] (0,17) grid (24,20);
\node [circle, draw, thick] (a1) at (1,19) {$0$};
\node [circle, draw, thick] (b1) at (5,19) {$1$};
\node [circle, draw, thick] (c1) at (9,19) {$2$};
\node [circle, draw, thick] (d1) at (13,19) {$3$};
\node [circle, draw, thick] (e1) at (17, 19){$4$};
\draw [-stealth, black, thick] (a1) -- node[below] {$0$} (b1);
\draw [-stealth, black, thick] (b1) -- node[below] {$0$} (c1);
\draw [-stealth, black, thick] (c1) -- node[below] {$1$} (d1);
\draw [-stealth, black, thick] (d1) -- node[below] {$0$} (e1);
\draw [-stealth, black, thick] (e1) -- node[below] {$1$} (20,19);
\draw [-stealth, black, thick] (b1) to [out=135, in=45] node[above] {$1$}(a1);
\draw [-stealth, black, thick] (e1) to [out=135, in=45] node[above] {$0$}(c1);
\path (c1) edge [loop above] node {$0$} (c1);

\end{tikzpicture}
\end{center}
\caption{An illustration of why $j\neq 2$ in the proof of Step 3.}
\label{case3graph}
\end{figure}
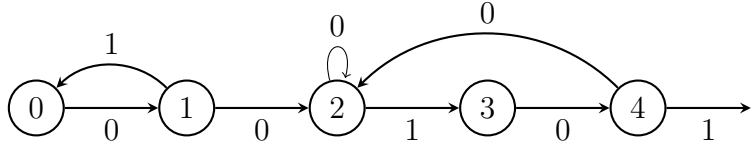

Thus, there are at most three possibilities: $j=0$, $j=1$ or $j=2$. 
We claim that $j=2$ is impossible. Indeed, if $j=2$, then $w_1w_2w_3w_4w_5=00101$, and the unlabeled version of the subgraph involving vertices $0,1,2$ is exactly $S_1$ (see Figure \ref{case3graph}), which contradicts the assumption we made at the beginning of this case.

Thus, we must have $j\in \{0,1\}$. 
Indeed, one directly verifies that $j=0$ corresponds to Case 3 of Figure \ref{all}, and $j=1$ corresponds to Case 6 of Figure \ref{all}, and the way to distinguish these two cases is to check the in-degree of vertices that can go to the vertex $I_2$: for Case 3, any vertex that  goes into the self-loop vertex must have in-degree $\geq 1$; for Case 6, the vertex $0$ goes into the self-loopvertex and, for a similar reason as in the proof of Step 2, the vertex $0$ must have in-degree $0$.

Having known how to distinguish Case 3 from Case 6, identifying the vertex labeled $0$ is straight forward: if we are in Case 3, then the vertex with a self-loop will be labeled $0$; if we are in Case 6, then there is a unique vertex with in-degree $0$ and goes into $I_1$, this vertex must be labeled $0$.

{\em{Proof of Step 4}}: If $S_1, S_2, S_3$ do not appear in $G$, then $S_4$ must appear. Now, by checking all the six cases in Figure \ref{all}, we conclude that we must be in Case 4. But in Case 4, $S_4$ must appear at vertices $0, 1$ and $2$. Thus, the vertex $I_1$ in $S_4$ must be labeled $0$.

%
%
%
%
%
%
%
%

\end{document}